\documentclass{amsart}
\usepackage{latexsym}
\usepackage{amsmath, amssymb}
\usepackage[all]{xy}

\theoremstyle{definition}
\newtheorem{dfs}{Definition}[section]

\newtheorem{rems}[dfs]{Remark}
\newtheorem{remss}[dfs]{Remarks}

\theoremstyle{plain}
\newtheorem{lms}[dfs]{Lemma}
\newtheorem{thms}[dfs]{Theorem}
\newtheorem{props}[dfs]{Proposition}

\newtheorem{cors}[dfs]{Corollary}

\newcommand{\vi}{\ensuremath{\mathcal{V}\mathrm{I}}}
\newcommand{\dr}{\mathrm{dr}\,}
\newcommand{\be}{{\bf1}}
\newcommand{\id}{\mathrm{id}}

\begin{document}

\title[The Elliott conjecture for Villadsen algebras]{The Elliott conjecture for Villadsen algebras \\ of the first type}
\author{Andrew S. Toms}
\address{Department of Mathematics and Statistics, York University, 
Toronto, Ontario, Canada, M3J 1P3}
\email{atoms@mathstat.yorku.ca}
\author{Wilhelm Winter}
\address{Mathematisches Institut der Universit\"at M\"unster\\
Einsteinstr. 62\\ D-48149 M\"unster}
\email{wwinter@math.uni-muenster.de}

\date{\today}

\thanks{{\it Research partly supported by:} Deutsche Forschungsgemeinschaft
(through the SFB 478), \\
\indent EU-Network Quantum Spaces - Noncommutative Geometry
(Contract No.\ HPRN-CT-2002-\indent00280), and an NSERC Discovery Grant}

\begin{abstract}
We study the class of simple C$^*$-algebras introduced by Villadsen in his pioneering work
on perforated ordered $\mathrm{K}$-theory.  We establish six equivalent characterisations of the proper subclass 
which satisfies the strong form of Elliott's classification conjecture:  two $C^*$-algebraic ($\mathcal{Z}$-stability and approximate
divisibility), 
one $\mathrm{K}$-theoretic (strict comparison of positive elements), and three topological (finite decomposition rank,
slow dimension growth, and bounded dimension growth).
The equivalence of $\mathcal{Z}$-stability and strict comparison constitutes a stably 
finite version of Kirchberg's characterisation of purely infinite $C^*$-algebras.  The other equivalences
confirm, for Villadsen's algebras, heretofore conjectural relationships between various notions of 
good behaviour for nuclear C$^*$-algebras.  
\end{abstract}

\maketitle

\section{Introduction}
The classification theory of norm-separable C$^*$-algebras began with Glimm's study
of UHF algebras in 1960 (\cite{Gli:UHF}), and was expanded by Bratteli's 1972 classification of 
approximately finite-dimensional (AF) algebras via certain directed graphs (\cite{Bra:AF}).  It was
with the work of Elliott, however, that the theory grew exponentially.  His classification
of both AF and A$\mathbb{T}$ algebras of real rank zero via their scaled ordered $\mathrm{K}$-theory
suggested a deep truth about the structure of separable and nuclear C$^*$-algebras (\cite{Ell:AF}, \cite{Ell:AT}).  
He articulated this idea in the late 1980s, and formalised it in his 1994 ICM address:  simple, separable, and nuclear
C$^*$-algebras should be classified up to $*$-isomorphism by their topological $\mathrm{K}$-theory
and traces.  This prediction came to be known as the Elliott conjecture.  

The 1990s and early 2000s saw Elliott's conjecture confirmed in remarkable generality.  Kirchberg and Phillips
established it for purely infinite C$^*$-algebras satisfying the Universal Coefficient Theorem (\cite{Kir:pi}, \cite{Phi:pi}),
and Lin did the same for his C$^*$-algebras of tracial topological rank zero (\cite{Lin:tr0}).  Elliott,
Gong, and Li, confirmed the conjecture for unital approximately homogeneous (AH) algebras of bounded
dimension growth (\cite{EGL}).  These results cover many natural examples of C$^*$-algebras, including those arising
from certain graphs, dynamical systems, and shift spaces.

In the midst of these successes, Villadsen produced a strange thing:  a simple, separable, and nuclear
C$^*$-algebra whose ordered $\mathrm{K}_0$-group was perforated, i.e., contained a non-positive element
$x$ such that $nx$ was positive and non-zero for some $n \in \mathbb{N}$ (\cite{V1}).  (This answered a long-standing
question of Blackadar concerning the comparison theory of projections in C$^*$-algebras.)  
The techniques used by Villadsen to study the $\mathrm{K}$-theory of his algebra were drawn from differential topology, 
and it took time for the functional analysts of the classification community to digest them.  
Then, in 2002, R\o rdam found a way to adapt Villadsen's techniques
to construct the first counterexample to the Elliott conjecture (\cite{R3}).  
Other counterexamples followed (\cite{To1}, \cite{To2}).

The success of Elliott's conjecture, however, is no accident.  It is a deep and fascinating phenomenon, and  
one must ask whether there is a regularity property lurking in those algebras for which
the Elliott conjecture is confirmed.  Various candidates exist:  stability under tensoring
with the Jiang-Su algebra $\mathcal{Z}$, finite decomposition rank, and, for approximately
subhomogeneous (ASH) algebras, the notion of strict slow dimension growth.  The first property
--- known as $\mathcal{Z}$-stability --- is perhaps the most natural candidate, since tensoring
with $\mathcal{Z}$ does not affect $\mathrm{K}$-theory or traces of a simple unital C$^*$-algebra
with weakly unperforated $\mathrm{K}_0$-group.  Elliott's conjecture thus predicts that all
such algebras will be $\mathcal{Z}$-stable.  It is this very prediction which forms the basis
for the counterexamples of R\o rdam and the first named author:  one produces pairs of simple unital
C$^*$-algebras with weakly unperforated $\mathrm{K}_0$-groups, one of which is not $\mathcal{Z}$-stable.
These examples have legitimised the assumption of $\mathcal{Z}$-stability in Elliott's classification
program, leading to the wide-ranging classification theorem of the second named author
for ASH algebras of real rank zero (\cite{Wi6}).  

The problem with $\mathcal{Z}$-stability in relation to Elliott's classification program is that  
its ability to {\it characterise} those algebras which are amenable to classification is an article of faith.
In all cases where $\mathcal{Z}$-stability is sufficient for classification (e.g., simple unital ASH algebras 
of real rank zero), it may also be automatic;  when it is known to be necessary for classification 
(e.g., AH algebras), it is not known to suffice.
In this paper we prove that $\mathcal{Z}$-stability does characterise those algebras which satisfy Elliott's 
conjecture in an ambient class where the assumption of $\mathcal{Z}$-stability is truly necessary.  The class considered
is at once substantial and the natural starting point for establishing such a characterisation:  Villadsen's algebras.
In fact we will prove much more.  
$\mathcal{Z}$-stability is not only the hoped for necessary and sufficient condition for classification, 
but is furthermore equivalent to a topological condition (finite decomposition rank) and to a
$\mathrm{K}$-theoretic condition (strict comparison of positive elements).  These three conditions,
all of which make sense for an arbitrary nuclear C$^*$-algebra, are
equivalent to three further conditions which are to varying extents native to the class of algebras we
consider:  approximate divisibility, slow dimension growth, and bounded dimension growth.

Some comments on our characterisations are in order.  Nuclear C$^*$-algebras can be viewed
from several angles.  They are evidently analytic objects, but can be seen as ordered algebraic
objects through their $\mathrm{K}$-theory, or as topological objects via the decomposition rank of 
Kirchberg and the second named author.  Our main result says that from each of these viewpoints, 
there is a natural way to characterise those C$^*$-algebras which satisfy the Elliott conjecture. 
The equivalence of $\mathcal{Z}$-stability, approximate divisibility, finite decomposition rank, 
slow dimension growth, and bounded dimension growth is a satisfying confirmation of the 
expectations of experts.  The equivalence of these conditions
with strict comparison of positive elements, however, is unexpected and exciting for several 
reasons.   First, the very idea of there being a $\mathrm{K}$-theoretic 
characterisation of those algebras which will satisfy Elliott's conjecture is new.  Second,
it is a condition that can be verified for large classes of examples generally suspected to be amenable
to classification (\cite{R4}, \cite{To5}).  Third, and most remarkably, this equivalence is 
a stably finite version of Kirchberg's celebrated characterisation of purely infinite C$^*$-algebras.

Our paper is organised as follows:  in Section 2 we recall the definitions of the regularity properties
which appear in our main result;  Section 3 introduces Villadsen algebras of the first type, and states
our main result;  Sections 4, 5, 6, and 7 contain the proof of the main result;  Section 8 gives some 
examples of non-$\mathcal{Z}$-stable Villadsen algebras.

\vspace{2mm}
\noindent
{\it Acknowledgement.}  Part of the work on this paper was carried out while the second named author visited
the first at the University of New Brunswick (UNB).  Both authors thank UNB and the Atlantic Centre for Operator
Algebras for their support. 

%%%%%%%%%%%%%%%%%%%%%%%%%%%%%%%%%%%%%%

\section{Preliminaries and notation}\label{prelim}

\subsection{AH algebras and dimension growth}

Below we recall the concepts of (separable unital) AH algebras and their dimension growth. 

\begin{dfs}
\label{sdg-def}
A separable unital $C^{*}$-algebra $A$ is called approximately homogeneous, or AH, if it can be written as an inductive limit 
\[
A= \lim_{i \to \infty} (A_{i},\phi_{i})
\]
where each $A_{i}$ is a  $C^{*}$-algebra of the form  
\[
A_{i}= \bigoplus_{j=1}^{m_{i}} p_{i,j}(\mathrm{C}(X_{i,j}) \otimes M_{r_{i,j}})p_{i,j}
\]
for natural numbers $m_{i}$ and  $r_{i,j}$,  compact metrisable spaces  $X_{i,j}$ and projections $p_{i,j} \in \mathrm{C}(X_{i,j}) \otimes M_{r_{i,j}}$. We refer to the inductive system $(A_{i},\phi_{i})_{i \in \mathbb{N}}$ as an AH decomposition for $A$. 

We say the AH decomposition $(A_{i},\phi_{i})_{\mathbf{N}}$ has slow dimension growth, if 
\[
\lim_{i \to \infty} \max_{j=1, \ldots,m_{i}} \frac{\dim X_{i,j}}{\mathrm{rank} \, p_{i,j}} = 0;
\]
it has very slow dimension growth, if  
\[
\lim_{i \to \infty} \max_{j=1, \ldots,m_{i}} \frac{(\dim X_{i,j})^{3}}{\mathrm{rank} \, p_{i,j}} = 0
\]
and it has bounded dimension growth, if 
\[
\sup_{i\in \mathbf{N}} \max_{j=1, \ldots,m_{i}} \dim X_{i,j} = d < \infty.
\]
The AH algebra $A$ has slow (very slow or bounded, respectively) dimension growth, if it has an AH decomposition which has slow (very slow or bounded, respectively) dimension growth.
\end{dfs}

\begin{rems} 
Slow dimension growth is obviously entailed by very slow dimension growth. Moreover, it is 
easy to see that if $A$ is simple, then bounded dimension growth implies very slow dimension 
growth. One of the remarkable results of \cite{G} says that, for simple AH algebras, very 
slow dimension growth also implies bounded dimension growth.  
\end{rems}

\subsection{Approximate divisibility and the Jiang--Su algebra}
Let $p$, $q$ and $n$ be natural numbers with $p$ and $q$ dividing $n$.  
$C^*$-algebras of the form
\begin{displaymath}
I[p,n,q] = \{f \in M_n(\mathrm{C}([0,1])) \, | \, f(0)=\be_{n/p} \otimes a, 
f(1)=b \otimes \be_{n/q},
a \in M_p, b \in M_q\}
\end{displaymath}
are commonly referred to as dimension drop intervals.  If $n=pq$ and $\mathrm{gcd}(p,q) 
= 1$, then the dimension drop
interval is said to be prime.  

In \cite{JS1}, Jiang and Su construct a $C^{*}$-algebra $\mathcal{Z}$, which is the unique simple unital inductive
limit of dimension drop intervals having $\mathrm{K}_0 = \mathcal{Z}$, $\mathrm{K}_1 
= 0$ and a unique normalised
trace.  It is a limit of prime dimension drop intervals where the matrix 
dimensions tend to infinity,
and there is a unital embedding of any prime dimension drop interval into 
$\mathcal{Z}$. Jiang and Su show that $\mathcal{Z}$ is strongly self-absorbing in the sense of \cite{TW1}. 

A $C^{*}$-algebra $A$ is said to be $\mathcal{Z}$-stable, if it is isomorphic to $A \otimes \mathcal{Z}$ (since $\mathcal{Z}$ is nuclear, there is no need to specify which tensor product we use). It was shown in \cite{JS1} and \cite{TW2} that all classes of simple $C^{*}$-algebras for which the Elliott conjecture has been verified so far consist of $\mathcal{Z}$-stable $C^{*}$-algebras. 

Using semiprojectivity of prime dimension drop intervals, it is not too hard to see that a separable unital $C^{*}$-algebra $A$ is $\mathcal{Z}$-stable if and only if the following holds (cf.\ \cite{TW2}): for any $n \in \mathbb{N}$ there is a sequence of unital completely positive contractions $\phi_{i}: M_{n} \oplus M_{n+1} \to A$ such that the restrictions of $\phi_{i}$ to $M_{n}$ and $M_{n+1}$ both preserve orthogonality (i.e., have order zero in the sense of \cite{Wi1}) and such that $\|[\phi_{i}((x,0)),\phi_{i}((0,y))]\| \to 0$ and $\|[\phi_{i}((x,y)),a]\| \to 0$ as $i $ goes to infinity for every $x \in M_{n}$, $y \in M_{n+1}$ and $a \in A$. This characterisation shows that $\mathcal{Z}$-stability generalises the concept of approximate divisibility:

Following \cite{BKR}, we say a separable unital $C^{*}$-algebra $A$ is approximately divisible, if for any $n \in \mathbb{N}$ there is a sequence of unital $*$-homomorphisms $\phi_{i}: M_{n} \oplus M_{n+1} \to A$ such that $\|[\phi_{i}(x),a]\| \to 0$ as $i $ goes to infinity for every $x \in M_{n} \oplus M_{n+1}$ and $a \in A$.

It was shown in \cite{TW2} that approximate divisibility indeed implies $\mathcal{Z}$-stability. The converse cannot hold in general (approximate divisibility asks for the existence of an abumdance of projections). However, using the classification result of \cite{EGL}, it was shown in \cite{EGL2}  that simple AH algebras of bounded dimension growth are approximately divisible.

\subsection{The Cuntz semigroup}
Let $A$ be a $C^*$-algebra, and let $\mathrm{M}_n(A)$ denote the $n \times n$ 
matrices whose entries are elements of $A$.  If $A = \mathbb{C}$, then we simply write $\mathrm{M}_n$.
Let $\mathrm{M}_{\infty}(A)$ denote the algebraic limit of the
direct system $(\mathrm{M}_n(A),\phi_n)$, where $\phi_n:\mathrm{M}_n(A) \to \mathrm{M}_{n+1}(A)$
is given by
\[
a \mapsto \left( \begin{array}{cc} a & 0 \\ 0 & 0 \end{array} \right).
\]
Let $\mathrm{M}_{\infty}(A)_+$ (resp. $\mathrm{M}_n(A)_+$)
denote the positive elements in $\mathrm{M}_{\infty}(A)$ (resp. $\mathrm{M}_n(A)$). 
Given $a,b \in \mathrm{M}_{\infty}(A)_+$, we say that $a$ is {\it Cuntz subequivalent} to
$b$ (written $a \precsim b$) if there is a sequence $(v_n)_{n=1}^{\infty}$ of
elements of $\mathrm{M}_{\infty}(A)$ such that
\[
||v_nbv_n^*-a|| \stackrel{n \to \infty}{\longrightarrow} 0.
\]
We say that $a$ and $b$ are {\it Cuntz equivalent} (written $a \sim b$) if
$a \precsim b$ and $b \precsim a$.  This relation is an equivalence relation,
and we write $\langle a \rangle$ for the equivalence class of $a$.  The set
\[
W(A) := \mathrm{M}_{\infty}(A)_+/ \sim
\] 
becomes a positively ordered Abelian monoid when equipped with the operation
\[
\langle a \rangle + \langle b \rangle = \langle a \oplus b \rangle
\]
and the partial order
\[
\langle a \rangle \leq \langle b \rangle \Leftrightarrow a \precsim b.
\]
In the sequel, we refer to this object as the {\it Cuntz semigroup} of $A$.

Given $a \in \mathrm{M}_{\infty}(A)_+$ and $\epsilon > 0$, we denote by 
$(a-\epsilon)_+$ the element of $C^*(a)$ corresponding (via the functional
calculus) to the function
\[
f(t) = \mathrm{max}\{0,t-\epsilon\}, \ t \in \sigma(a).
\]
(Here $\sigma(a)$ denotes the spectrum of $a$.)  

%The proposition below collects some facts about Cuntz subequivalence due 
%to Kirchberg and R{\o}rdam.

%\begin{props}[Kirchberg-R{\o}rdam (\cite{KR}), R{\o}rdam (\cite{R4})]\label{basics}
%Let $A$ be a $C^*$-algebra, and $a,b \in A_+$.  
%\begin{enumerate}
%\item[(i)] $(a-\epsilon)_+ \precsim a$ for every $\epsilon > 0$.
%\item[(ii)] The following are equivalent:
%\begin{enumerate}
%\item[(a)] $a \precsim b$;
%\item[(b)] for all $\epsilon > 0$, $(a-\epsilon)_+ \precsim b$;
%\item[(c)] for all $\epsilon > 0$, there exists $\delta > 0$ such that $(a-\epsilon)_+ \precsim (b-\delta)_+$.
%\end{enumerate}
%\item[(iii)] If $\epsilon>0$ and $||a-b||<\epsilon$, then $(a-\epsilon)_+ \precsim b$.
%\end{enumerate}
%\end{props} 

\subsection{Dimension functions and strict comparison}
Now suppose that $A$ is unital and stably finite, and denote by $\mathrm{QT}(A)$
the space of normalised 2-quasitraces on $A$ (v. \cite[Definition II.1.1]{BH}).
Let $S(W(A))$ denote the set of additive and order preserving maps $d$ from $W(A)$ to $\mathbb{R}^+$
having the property that $d(\langle 1_A \rangle) = 1$.
Such maps are called {\it states}.  Given $\tau \in \mathrm{QT}(A)$, one may 
define a map $d_{\tau}:\mathrm{M}_{\infty}(A)_+ \to \mathbb{R}^+$ by
\begin{equation}\label{ldf}
d_{\tau}(a) = \lim_{n \to \infty} \tau(a^{1/n}).
\end{equation}
This map is lower semicontinous, and depends only on the Cuntz equivalence class
of $a$.  It moreover has the following properties:
\vspace{2mm}
\begin{enumerate}
\item[(i)] if $a \precsim b$, then $d_{\tau}(a) \leq d_{\tau}(b)$;
\item[(ii)] if $a$ and $b$ are orthogonal, then $d_{\tau}(a+b) = d_{\tau}(a)+d_{\tau}(b)$;
\item[(iii)] $d_{\tau}((a-\epsilon)_+) \nearrow d_{\tau}(a)$ as $\epsilon \to 0$. 
\end{enumerate}
\vspace{2mm}
Thus, $d_{\tau}$ defines a state on $W(A)$.
Such states are called {\it lower semicontinuous dimension functions}, and the set of them 
is denoted $\mathrm{LDF}(A)$.  If $A$ has the property that $a \precsim b$ whenever $d(a) < d(b)$
for every $d \in \mathrm{LDF}(A)$, then we say that $A$ has 
{\it strict comparison of positive elements} or simply {\it strict comparison}.

\subsection{Decomposition rank}
Based on the completely positive approximation property for nuclear $C^{*}$-algebras, one may define a noncommutative version of covering dimension as follows:

\begin{dfs}[\cite{KW}, Definitions 2.2 and 3.1]
\label{dr-def}
Let $A$ be a separable $C^*$-algebra. 

(i) A completely positive map $\varphi : \bigoplus_{i=1}^s M_{r_i} \to A$ is $n$-decomposable, if there is a decomposition $\{1, \ldots, s\} = \coprod_{j=0}^n I_j$ such that the restriction of $\varphi$ to $\bigoplus_{i \in I_j} M_{r_i}$ preserves orthogonality for each $j \in \{0, \ldots, n\}$.

(ii) $A$ has decomposition rank $n$, $\dr A = n$, if $n$ is the least integer such that the following holds: Given $\{b_1, \ldots, b_m\} \subset A$ and $\epsilon > 0$, there is a completely positive approximation $(F, \psi, \varphi)$ for $b_1, \ldots, b_m$ within $\epsilon$ (i.e., $\psi:A \to F$ and $\varphi:F \to A$ are completely positive contractions and $\|\varphi \psi (b_i) - b_i\| < \epsilon$) such that $\varphi$ is $n$-decomposable. If no such $n$ exists, we write $\dr A = \infty$.  
\end{dfs}

This notion has good permanence properties; for example, it behaves well with respect to quotients, inductive limits, hereditary subalgebras, unitization and stabilization. It generalises topological covering dimension, i.e., if $X$ is a locally compact second countable space, then $\dr \mathrm{C}_0(X) = \dim X$; see \cite{KW} for details. Moreover, if $A$ is an AH algebra of bounded dimension growth, then $\dr A$ is finite.

%%%%%%%%%%%%%%%%%%%%%%%%%%%%%%%%%%%%%

\section{$\vi$ Algebras and the main result}

\subsection{Villadsen algebras of the first type}

The class of algebras we consider is an interpolated family of AH algebras.
At their simplest they are the UHF algebras of Glimm, while at their most complex they are the algebras
introduced by Villadsen in his work on perforated ordered $\mathrm{K}$-theory.   In between these extremes
they span the full spectrum of complexity for simple, separable, nuclear, and stably finite C$^*$-algebras.   
We call these algebras {\it Villadsen algebras of the first type} as they are defined by a generalisation of 
Villadsen's construction in \cite{V1}.  (Villadsen used a second and quite distinct construction
in his subsequent work on stable rank, cf.\ \cite{V2}.) 

Let $X$ and $Y$ be compact Hausdorff spaces.  Recall that a $*$-homomorphism 
\[
\phi:\mathrm{C}(X) \to \mathrm{M}_n \otimes \mathrm{C}(Y) 
\]
is said to be {\it diagonal} if it has the form
\[
f \mapsto \mathrm{diag}(f \circ \lambda_1,\ldots,f \circ \lambda_n),
\]
where $\lambda_i:Y \to X$ is a continuous map for each $1 \leq i \leq n$.  The maps $\lambda_1,\ldots,\lambda_n$
are called the {\it eigenvalue maps} of $\phi$.  Amplifications of diagonal maps are
also called diagonal.

\begin{dfs}\label{vmap}
Let $X$ be a compact Hausdorff space and $n,m,k \in \mathbb{N}$.  A unital diagonal $*$-homomorphism
\[
\phi: \mathrm{M}_n \otimes \mathrm{C}(X) \to \mathrm{M}_k \otimes \mathrm{C}(X^{\times m})
\]
is said to be a Villadsen map of the first type (a $\vi$ map) if each eigenvalue map is either a co-ordinate
projection or has range equal to a point.
\end{dfs}

\begin{dfs}\label{vonealg}
Let $X$ be a compact Hausdorff space, and let $(n_i)_{i=1}^{\infty}$ and $(m_i)_{i=1}^{\infty}$ be sequences
of natural numbers with $n_{1}=1$.  Fix a compact Hausdorff space $X$, and put $X_i = X^{\times n_i}$. 
A unital $C^{*}$-algebra $A$ is said to be a Villadsen algebra of the first type (a $\vi$ algebra), if it can be written as an inductive limit algebra 
\[
A \cong \lim_{i \to \infty}\left( \mathrm{M}_{m_i} \otimes \mathrm{C}(X^{\times n_i}), \phi_i \right)
\]
where each $\phi_i$ is a $\vi$ map. 
\end{dfs}
\noindent
We will refer to the inductive system in Definition \ref{vonealg} as
a {\it standard decomposition} for $A$ with {\it seed space} $X_1 (=X)$.
Clearly, such decompositions are not unique.

For $j>i$, put
\[
\phi_{i,j} = \phi_{j-1} \circ \cdots \circ \phi_i.
\]
Let $N_{i,j}$ be the number of distinct co-ordinate projections
from $X_{j} = X_i^{\times n_{j}/n_i}$ to $X_i$ occuring as eigenvalue
maps of $\phi_{i,j}$, and let $M_{i,j}$ denote the multiplicity (number of 
eigenvalue maps) of $\phi_{i,j}$.  Notice that 
\[
M_{i,j} = M_{j-1,j} M_{i,j-1}, \ \mathrm{that} \ N_{i,j} = N_{j-1,j} N_{i,j-1}
\]
and that
\[
0 \le \frac{N_{i,j}}{M_{i,j}} \le 1.
\]
From these relations it follows in particular that the sequence
\[
\left(  \frac{N_{i,j}}{M_{i,j}} \right)_{j>i}
\]
is decreasing and converges for any fixed $i$.   

$A$ is said to have {\it  slow (very slow, or bounded) dimension growth as a $\vi$ algebra}, if it admits a 
standard decomposition as above which has slow (very slow, or bounded, respectively) dimension growth in the sense of Definition \ref{sdg-def}.

\begin{remss} {\rm Despite its simple definition, the class of $\vi$ algebras is surprisingly broad:\\

\noindent
$\bullet$ By taking $X_1 = \{*\}$, we can obtain any UHF algebra.
If instead we take $X_1$ to be a finite set, then we obtain a good supply of AF algebras.\\

\noindent
$\bullet$ With each $X_i$ equal to a disjoint union of finitely many circles, we obtain a large collection
of A$\mathbb{T}$ algebras of real rank zero and real rank one.\\

\noindent
$\bullet$ If each $X_i$ is equal to the same compact Hausdorff space $X$, then we obtain the class of Goodearl algebras.\\

\noindent
$\bullet$ If we impose the condition that $n_i/m_i \to 0$, then we obtain AH algebras of slow dimension growth
exhibiting a full range of complexity in their Elliott invariants:  torsion in $\mathrm{K}_0$ or 
$\mathrm{K}_1$, and arbitrary tracial state spaces.  \\

\noindent
$\bullet$ By taking ``most'' of the eigenvalue maps in each $\phi_i$ to be distinct co-ordinate projections 
and setting $X_1 = \mathrm{S}^2$
we obtain Villadsen's example of a simple, separable, and nuclear $C^*$-algebra with perforated ordered
$\mathrm{K}_0$-group (\cite{V1}). A variation on Villadsen's construction yields the counterexample to 
Elliott's classification conjecture discovered by the first named author in \cite{To2}. 
}
\end{remss}

\vspace{2mm}
The first three examples above are special cases of the fourth, and the latter is a class of
algebras for which the Elliott conjecture can be shown to hold.  Proving this, however, requires
both the most powerful available classification results for stably finite C$^*$-algebras and 
the detailed analysis of $\vi$ algebras provided in the sequel.  Thus, 
from the standpoint of trying to confirm the Elliott conjecture, $\vi$ algebras
are no less complex than the class of all simple unital AH algebras.  The fifth example
demonstrates that $\vi$ algebras include non-$\mathcal{Z}$-stable algebras which, in general,
cannot be detected with classical $\mathrm{K}$-theory.  The sequel will show that there are in fact a 
tremendous number of such $\vi$ algebras (see Section 8).

\subsection{The main theorem}

\begin{thms}\label{main}
Let $A$ be a simple $\vi$ algebra admitting a standard decomposition with 
seed space a finite-dimensional CW complex.  The following are equivalent:
\begin{enumerate}
\item[(i)] $A$ is $\mathcal{Z}$-stable;
\item[(ii)] $A$ has strict comparison of positive elements;
\item[(iii)] $A$ has finite decomposition rank;
\item[(iv)] $A$ has slow dimension growth (as an AH algebra);
\item[(v)] $A$ has bounded dimension growth (as an AH algebra);
\item[(vi)] $A$ is approximately divisible.
\end{enumerate}
If, moreover, $A$ has real rank zero, then $A$ satisfies the equivalent conditions
above. 
\end{thms}

For Theorem \ref {main}, the following implications are already known:
\[
\mathrm{(v)} \Rightarrow \mathrm{(vi)} \Rightarrow \mathrm{(i)} \Rightarrow \mathrm{(ii)}; \ \ 
\mathrm{(v)} \Rightarrow \mathrm{(iii)}; \ \ \mathrm{(iv)} \Rightarrow \mathrm{(ii)}.
\]
More precisely, $\mathrm{(v)} \Rightarrow \mathrm{(vi)}$ was shown in \cite{EGL2} (based on 
the results of \cite{EGL}), $\mathrm{(vi)} \Rightarrow \mathrm{(i)}$ is a result of \cite{TW2}, 
$\mathrm{(i)} \Rightarrow \mathrm{(ii)}$ is \cite[Corollary~4.6]{R4}, $\mathrm{(iv)} \Rightarrow
\mathrm{(ii)}$ is \cite[Corollary 4.6]{To5}, and $\mathrm{(v)} \Rightarrow 
\mathrm{(iii)}$ is an easy observation of \cite{KW}. We will prove $\mathrm{(ii)} \Rightarrow 
\mathrm{(iv)}$, $\mathrm{(ii)} \Rightarrow \mathrm{(v)}$, $\mathrm{(iii)} \Rightarrow \mathrm{(ii)}$, 
and the statement about real rank zero.

\begin{remss} $\bullet$ In the absence of Theorem \ref{main}, no two of conditions $\mathrm{(i)-(vi)}$ are 
known to be equivalent for an arbitrary simple unital AH algebra. However, $\mathrm{(i)}$, 
$\mathrm{(ii)}$, $\mathrm{(iv)}$, $\mathrm{(v)}$ and $\mathrm{(vi)}$ are known to be equivalent for simple 
unital AH algebras of real rank zero (a combination of results from \cite{EGL2}, \cite{TW2}, 
\cite{Wi6}, \cite{Lin:tr0} and \cite{El4}).  One of the central open questions in the classification
theory of nuclear C$^*$-algebras is whether any of these conditions is actually necessary in
the real rank zero case.  Theorem \ref{main} makes some progress on this question --- see the third remark below. \\  

\noindent
$\bullet$ Conditions $\mathrm{(i)-(iii)}$ should remain equivalent in much larger classes of simple, 
separable, nuclear, and stably finite $C^*$-algebras. Conditions $\mathrm{(iv)}$, $\mathrm{(v)}$,
and $\mathrm{(vi)}$ cannot be expected to hold in general (conditions $\mathrm{(iv)}$ 
and $\mathrm{(v)}$ exclude  non-AH algebras, and 
$\mathrm{(vi)}$ excludes  projectionless algebras, such as $\mathcal{Z}$ itself), but it remains 
possible that $\mathrm{(i)-(vi)}$ are equivalent for simple unital ASH algebras, when $\mathrm{(iv)}$ and $\mathrm{(v)}$  
are adapted in the obvious way to this class. \\

\noindent
$\bullet$  It is not known whether real rank zero implies conditions $\mathrm{(i)-(vi)}$ above for larger 
classes of simple nuclear $C^{*}$-algebras (clearly, $\mathrm{(iii)}$, $\mathrm{(iv)}$, and $\mathrm{(v)}$ can only 
hold in the stably finite case). Every existing classification theorem for real rank zero C$^*$-algebras assumes at least one
of conditions $\mathrm{(i)-(vi)}$.  It thus remarkable that in the class of  $\vi$ algebras, 
real rank zero entails classifiability without assuming any of these conditions.  \\

\noindent
$\bullet$ Our proof of $\mathrm{(iv)} \Rightarrow \mathrm{(v)}$ is the first instance of such among simple
unital AH algebras of unconstrained real rank. \\

\noindent
$\bullet$ The proof of Theorem \ref{main} yields new examples of simple C$^*$-algebras with infinite decomposition
rank.  Previous examples all had a unique trace, while our examples can exhibit a wide variety of structure
in the tracial state space. (See Section 8.)\\

\noindent
$\bullet$ Simple $\vi$ algebras all have stable rank one by an argument similar to that of 
\cite[Proposition 10]{V1}.  They may, however, have quite fast dimension growth ---
\cite[Theorem 5.1]{To6} exhibits a simple $\vi$ algebra for which every AH decomposition has 
the property that
\[
\liminf_{i \to \infty} \max_{j=1, \ldots,m_{i}} \frac{\dim X_{i,j}}{\mathrm{rank} \, p_{i,j}} = \infty.
\]
\end{remss}

\vspace{2mm}
For completeness, we note that the class of $\vi$ algebras which satisfy the conditions of 
Theorem \ref{main} indeed satisfies the Elliott conjecture. Let $\mathrm{Ell}(\bullet)$ 
denote the Elliott invariant of a unital, exact, and stably
finite $C^*$-algebra. We then have:

\begin{cors}[Gong (\cite{G}), Elliott-Gong-Li (\cite{EGL})]
Let $A$ and $B$ be simple $\vi$ algebras as in Theorem \ref{main} which satisfy 
conditions $\mathrm{(i)-(vi)}$, and suppose that there is an isomorphism
\[
\phi:\mathrm{Ell}(A) \to \mathrm{Ell}(B).
\]
Then, there is a $*$-isomorphism $\Phi:A \to B$ which induces $\phi$.
\end{cors}

\subsection{An analogue of Kirchberg's First Geneva Theorem}
The most interesting aspect of Theorem \ref{main} is that it provides an analogue among $\vi$ algebras
of Kirchberg's characterisation of purely infinite algebras.  The latter states that for a simple,
separable, and nuclear C$^*$-algebra $A$ we have
\[
A \otimes \mathcal{O}_{\infty} \cong A \Leftrightarrow \mathrm{A \ is \ purely \ infinite}.
\]
If we suppose that $A$ is {\it a priori} traceless, then a result of R\o rdam (see \cite{R4})
says that $\mathcal{Z}$-stability and $\mathcal{O}_{\infty}$-stability are equivalent, and
the definition of strict comparison reduces to the very definition of pure infiniteness.
Thus, we see that Kirchberg's characterisation is equivalent to the statement
\[
A \otimes \mathcal{Z} \cong A \Leftrightarrow \mathrm{A \ has \ strict \ comparison \ of \ positive \ elements}.
\]
This statement makes sense even if $A$ has a trace, and is moreover true for the simple $\vi$ algebras
of Theorem \ref{main}.  Were the statement to hold for all 
simple, separable, and nuclear C$^*$-algebras --- a distinct possibility --- it would be a deep
and striking generalisation of Kirchberg's characterisation.  In light of this possibility,
we suggest that simple, finite, and $\mathcal{Z}$-stable C$^*$-algebras be termed ``purely finite''.

%%%%%%%%%%%%%%%%%%%%%%%%%%%%%%%%%%%%%%%%%%%

%%%%%%%%%%%%%%%%%%%%%%%%%%%%%%%%%%%%%%%%

\section{Villadsen's obstruction in the Cuntz semigroup}\label{villob}

In this section we prove that under a technical assumption, a simple $\vi$ algebra
fails to have strict comparison of positive elements.  We shall see later that this failure
is dramatic enough to ensure that the algebra also has infinite decomposition rank.

\subsection{Vector bundles and characteristic class obstructions}
All vector bundles considered in this paper are topological and complex.
For any connected topological space $X$, we let $\theta_l$ denote the
trivial vector bundle of fibre dimension $l \in \mathbb{N}$. 
If $\omega$ is a vector bundle over $X$, then
we denote by $\oplus_{i=1}^k \omega$ the $k$-fold Whitney sum of $\omega$ 
with itself, and by $\omega^{\otimes k}$ its $k$-fold external tensor product
(over $X^k$).  We use $\mathrm{rank}(\omega)$ to denote 
the fibre dimension of $\omega$.  If $Y$ is a second topological space and 
$f:Y \to X$ is continuous, then $f^*(\omega)$ denotes the induced bundle over $Y$. 
By Swan's Theorem, $\omega$ can be represented by a (non-unique) projection in
a matrix algebra over $\mathrm{C}(X)$;  we will use vector bundles and projections
interchangeably in the sequel.

Recall that the Chern class $c(\omega)$ is an element
of the integral cohomology ring $H^*(X)$ of the form
\[
c(\omega) = \sum_{i=0}^{\infty} c_i(\omega),
\]
where $c_i(\omega) \in H^{2i}(X)$ and $c_i(\omega) = 0$ whenever $i > \mathrm{rank}(\omega)$.  
Let $\gamma$ be a second vector bundle over
$X$.  We will make use of the following properties of the Chern class:
\begin{enumerate}
\item[(i)] $c(\theta_l) = 1 \in H^0(X)$; 
\item[(ii)] $c(\gamma \oplus \omega) = c(\gamma)c(\omega)$, where the product is the cup product; 
\item[(iii)] If $Y$ is another topological space and $f:Y \to X$ is continuous, then
$c(f^*(\omega)) = f^*(c(\omega))$.
\end{enumerate}

Let $\xi$ be the Hopf line bundle over $\mathrm{S}^2$.
The following Chern class obstruction argument, due essentially to Villadsen, shows
that $\theta_{k}$ is not isomorphic to a sub-bundle of $\oplus_{i=1}^l \xi^{\otimes l}$
whenever $1 \leq k < l$.
The top Chern class $c_l(\oplus_{i=1}^l \xi^{\otimes l})$ (equal, in this
case, to the Euler class of $\oplus_{i=1}^l \xi^{\otimes l}$) is not zero by \cite[Proposition 3.2]{R3}.
If $\theta_{k}$ is isomorphic to a sub-bundle of $\oplus_{i=1}^l \xi^{\otimes l}$, then there
exists a vector bundle $\gamma$ of rank $l-k$ over $(\mathrm{S}^2)^l$
such that
\[
\theta_{k} \oplus \gamma \cong \oplus_{i=1}^l \xi^{\otimes l}.
\]
Applying the Chern class to this equation yields
\[
c(\gamma) = c(\oplus_{i=1}^l \xi^{\otimes l}).
\]
But then $c_l(\gamma) \neq 0$, 
contradicting the fact that $c_i(\gamma) = 0$ whenever $i > \mathrm{rank}(\gamma) = k$.

We review for future use some structural aspects of the integral cohomology ring 
$H^*((\mathrm{S}^2)^{n})$.  It is well known that 
\[
H^0(\mathrm{S}^2) \cong H^2(\mathrm{S}^2) \cong \mathbb{Z}
\]
and
\[
H^i(\mathrm{S}^2) = 0, \ i \neq 0,2.
\]
It follows from the K{\"u}nneth formula that 
\[
H^*((\mathrm{S}^2)^{ n}) \cong H^*(\mathrm{S}^2)^{\otimes n}
\]
as graded rings.  Let $e_i$ denote the generator of $H^2(\mathrm{S^2})$
in the $i^{\mathrm{th}}$ tensor factor of $H^*(\mathrm{S}^2)^{\otimes n}$.  Then,
\[
H^*((\mathrm{S}^2)^{ n}) \cong \langle 1,e_1,\ldots,e_n \rangle/R.
\]
where 
\[
R = \{e_i^2=0 \ | \ 1 \leq i \leq n\}.
\]
If $n = Nl$ for some $N \in \mathbb{N}$, then
\[
H^*\left((\mathrm{S}^2)^{ Nl}\right) = H^*\left((\mathrm{S}^2)^{ l}\right)^{\otimes N}.
\]
Let $e_{i,j}$ denote the generator of the $i^{\mathrm{th}}$ copy of 
$H^2(\mathrm{S}^2)$, $i \in \{1,\ldots,l\}$, in the $j^{\mathrm{th}}$
tensor factor of the right hand side above. 

\subsection{A failure of strict comparison in $\mathrm{C}(X)$}
Villadsen's Chern class obstruction argument may be viewed as a statement
about projections in a matrix algebra over $\mathrm{C}(X)$.  We present
below an analogue of his argument for certain non-projections 
in $\mathrm{M}_n(\mathrm{C}(X))$. 

Let $X$ be a CW-complex with $\mathrm{dim}(X) \geq 6$,
and let there be given a natural number $l$ satisfying $2 \leq l \leq \lfloor \mathrm{dim}(X)/3 \rfloor$.
Choose an open set $O \subseteq X$ homeomorphic to $(-1,1)^{\mathrm{dim}(X)}$.
Define
\[
\tilde{A} := \{ x \in (-1,1)^3 \ | \ \mathrm{dist}( x,(0,0,0)) = 1/2 \} \cong \mathrm{S}^2
\]
and
\[
\tilde{B} := \{ x \in (-1,1)^3 \ | \ 1/3 < \mathrm{dist}( x,(0,0,0)) < 2/3 \},
\]
and let $\pi:\tilde{B} \to \tilde{A}$ be the continuous projection along rays emanating from 
$(0,0,0)$.  Now define a closed subset
\[
A = \tilde{A}^l \times \{0\}^{\mathrm{dim}(X)-3l}
\]
and an open subset
\[
B = \tilde{B}^l \times (-1,1)^{\mathrm{dim}(X)-3l}
\]
of $O$.  Define a continuous map $\Pi:B \to A$ by 
\[
\Pi = \underbrace{\pi \times \cdots \times \pi}_{l \ \mathrm{times}} \times
\underbrace{ ev_0 \times \cdots \times ev_0}_{\mathrm{dim}(X)-3l \ \mathrm{times}},
\]
where $ev_0(x) = 0$ for every $x \in (-1,1)$.   
Let $f:X \to [0,1]$ be a continuous map which is identically one
on $A$ and identically zero off $B$.

Notice that $A \cong (\mathrm{S}^2)^l$, so $\xi^{\otimes l}$ may be viewed
as a vector bundle over $A$.  Define positive elements 
\begin{equation}\label{P}
P = f \cdot \Pi^*(\oplus_{i=1}^l \xi^{\otimes l}) 
\end{equation}
and
\begin{equation}\label{T}
\Theta_k = f \cdot \theta_{k}
\end{equation}
in $\mathrm{M}_{\infty}(\mathrm{C}(X))$.
For every $x \in X$ and $n \in \mathbb{N}$ we have either
\[
\mathrm{rank}(P(x)) = \mathrm{rank}(\Theta_k(x)) = 0
\]
or
\[
\mathrm{rank}(P(x)) = l, \ \mathrm{rank}(\Theta_k(x)) = k.
\]
If $\tau \in \mathrm{T}(\mathrm{C}(X))$ and $a \in \mathrm{M}_{\infty}(\mathrm{C}(X))_+$,
then $d_{\tau}(a)$ is obtained by integrating the rank function of $a$ against the probability
measure on $X$ corresponding to $\tau$.  Thus, if $k < l$, we have
\[
d_{\tau}(\Theta_k) < d_{\tau}(P), \ \forall \tau \in \mathrm{T}(A).
\]

On the other hand, $\langle \Theta_k \rangle \nleq \langle P \rangle$.  To see this suppose, on the
contrary, that there exists a sequence $(v_i)_{i=1}^{\infty}$ in
$\mathrm{M}_{\infty}(\mathrm{C}(X))$ such that 
\[
||v_i P v_i^* - \Theta_k|| \stackrel{i \to \infty}{\longrightarrow} 0.
\]
Then, the same is true upon restriction to $A \subseteq X$,
i.e., 
\[
\Theta_k|_{A} \cong \theta_{k} \precsim \oplus_{i=1}^l \xi^{\otimes l} \cong P|_{A}.
\]
This amounts to saying that $\theta_k$ is isomorphic to a sub-bundle of
$\oplus_{i=1}^l \xi^{\otimes l}$, contradicting our choice of $\xi$.

We are now ready to prove a key lemma.  Its proof is inspired by the proof of 
\cite[Theorem 1.1]{To2}.

%\begin{lms}\label{ob}
%Let $A$ be a $\vi$ algebra with standard decomposition $(A_i,\phi_i)$.  Suppose 
%that there exists $i \in \mathbb{N}$ such that $\lfloor \mathrm{dim}(X_i)/3 \rfloor \geq 3$,
%and $j \geq i$ such that
%\[
%\frac{N_{i,j}}{M_{i,j}} \geq \frac{kl-1}{kl}
%\]
%for some natural number $l$ satisfying $3 \leq l \leq \lfloor \mathrm{dim}(X_i)/3 \rfloor$
%and some $2 \leq k < l$.
%Let $P$ and $\Theta_k$ be the positive elements in a matrix algebra over $\mathrm{C}(X_i)$ constructed 
%as above.  Then, for every $v \in \mathrm{M}_{\infty}(A_{j})$ we have
%\[
%|| v \phi_{i,j}(P) v^* - \phi_{i,j}(\Theta_k) || \geq 1/2.
%\]
%\end{lms}

\begin{lms}\label{ob}

Let $A$ be a simple $\vi$ algebra with standard decomposition $(A_i,\phi_i)$
and seed space a CW-complex $X_1$ of dimension strictly greater than zero.
Suppose that for any $\epsilon > 0$ there exists $i \in \mathbb{N}$
such that
\begin{equation}
\label{N-M-sequence-1}
\frac{N_{i,j}}{M_{i,j}} > 1-\epsilon, \ \forall j >i.
\end{equation}

Then, for any $n \in \mathbb{N}$ there exist pairwise orthogonal elements 
$a, b_1,\ldots, b_n \in \mathrm{M}_{\infty}(A)_+$
such that for each $s \in \{1,\ldots,n\}$
\[
d_{\tau}(a) < d_{\tau}(b_s), \ \forall \tau \in \mathrm{T}(A), 
\]
and
\[
\langle a \rangle \nleq \langle b_1 \rangle + \cdots + \langle b_n \rangle.
\]
In particular, $A$ does not have strict comparison of positive elements.
\end{lms}

\begin{proof}
First observe that the simplicity of $A$ combined with the nonzero dimension
of $X_1$ imply that $m_i \to \infty$ as $i \to \infty$ --- the number of point evaluations
appearing as eigenvalue maps in $\phi_{i,j}$ is unbounded as $j \to \infty$.  It then follows from 
our assumption on $N_{i,j}/M_{i,j}$ that $\mathrm{dim}(X_i) \to \infty$ as $i \to \infty$.
We may thus assume that $\mathrm{dim}(X_i) \neq 0$, $\forall i \in \mathbb{N}$.
Since $\phi_i$ always contains at least one eigenvalue map which is not a point evaluation,
it is injective.

Let $n \in \mathbb{N}$ be given.  Find, using the hypotheses of the lemma,
an $i \in \mathbb{N}$ such that
\[
\frac{N_{i,j}}{M_{i,j}} > \frac{6n-1}{6n}, \ \forall j >i.
\]
Since $N_{i,j}/M_{i,j}$ increases with increasing $i$, we may assume that
$i$ is large enough to permit the construction of the element 
\[
f \cdot \Pi^*(\oplus_{i=1}^3 \xi^{\otimes 3n})
\]
--- this is just the element $P$ of equation (\ref{P}), with the number of
direct summands altered.  This element will be
our $b_1$.  (The maps $\phi_i$ are injective, so we identify forward images
in the inductive sequence.)  For $b_2,\ldots,b_n$ we simply take mutually orthogonal copies of 
$b_1$.  Let $a$ be $\Theta_2$, chosen orthogonal to each of $b_1,\ldots,b_n$.
We have by construction that
\[
d_{\tau}(a) < d_{\tau}(b_s), \ \forall \tau \in \mathrm{T}(A_i), \ 1 \leq s \leq n.
\]
This statement then holds for each $\tau \in \mathrm{T}(A)$ since
\[
d_{\tau}(\phi_{i\infty}(a)) = d_{\phi_{i\infty}^{\sharp}(\tau)}(a).
\] 

It remains to prove that 
\[
\langle a \rangle \nleq \langle b_1 \rangle + \cdots + \langle b_n \rangle = \langle b_1 + \cdots + b_n \rangle.
\]
Notice that $b_1 + \cdots +b_n$, viewed as an element of $\mathrm{M}_{\infty}(A_i)$, is simply
the element $P$ of equation (\ref{P}) with parameter $l = 3n$.  Thus, with this choice of $l$,
we are in fact trying to prove that 
\[
\langle \phi_{i\infty}(\Theta_2) \rangle \nleq \langle \phi_{i\infty}(P) \rangle.
\]
It will suffice to prove that for each $j > i$ and $v \in \mathrm{M}_{\infty}(A_j)$
\[
|| v \phi_{i,j}(P) v^* - \phi_{i,j}(\Theta_2) || \geq 1/2.
\]

Let $S$ be the set of eigenvalue maps of $\phi_{i,j}$.  $S$ is the
disjoint union of the set $S_1$ of eigenvalue maps which are co-ordinate
projections and the set $S_2$ of eigenvalue maps which are point evaluations.
(The fact that $\mathrm{dim}(X_i) \neq 0$ ensures that 
$S_1 \cap S_2 = \emptyset$.)  Note that $|S_1| = N_{i,j}$.
For $\lambda \in S_1$, let $m(\lambda)$ denote the number of times 
that $\lambda$ occurs as an eigenvalue map of some $\phi_{i,j}$.

Write $\phi_{i,j} = \gamma_1 \oplus \gamma_2$, where $\gamma_1$ is a $\vi$
map corresponding to the eigenvalue maps of $\phi_{i,j}$ contained in $S_1$,
and $\gamma_2$ corresponds similarly to the elements of $S_2$.
By construction, $\gamma_2(P)$ is a constant positive matrix-valued function
over $X_{j+1}$.  Put $\tilde{P} = \gamma_1(\mathbf{1}_{\mathrm{M}_n(\mathrm{C}(X_i))}) \oplus \gamma_2(P)^{1/2}$,
and $q = \lim_{n \to \infty} \gamma_2(P)^{1/n}$.
It follows that
\[
\phi_{i,j}(P) = \gamma_1(P) \oplus \gamma_2(P) = \tilde{P} (\gamma_1(a) \oplus q) \tilde{P},
\]
and that the projection $q$ corresponds to a trivial vector bundle.

Suppose that there exists $v \in \mathrm{M}_{\infty}(A_{j})$ such that
\[
|| v \phi_{i,j}(P) v^* - \phi_{i,j}(\Theta_2) || < 1/2.
\] 
Then, 
\[
|| v \tilde{P} (\gamma_1(P) \oplus q) \tilde{P} v^* - (\gamma_1(\Theta_2) \oplus \gamma_2(\Theta_2)) || < 1/2.
\]
Cutting down by $\gamma_1(\mathbf{1}_{A_i})$ and setting $w = \gamma_1(\mathbf{1}_{A_i}) v \tilde{P}$
we have
\begin{equation}\label{est1}
|| w (\gamma_1(P) \oplus q) w^* - \gamma_1(\Theta_2)|| < 1/2,
\end{equation}
and this estimate holds {\it a fortiori} over any closed subset of $X_{j}$.

Fix a point $x_0 \in X_i$ and let $C$ be the closed subset of $X_{j}=X_i^{\times n_{j}/n_i}$ 
consisting of those $(n_{j}/n_i)$-tuples which are equal to $x_0$ in those co-ordinates 
which are not the range of an element of $S_1$, and whose remaining co-ordinates belong
to $A \subseteq X_i$.  Notice that 
\[
C \cong A^{\times l N_{i,j}} \cong (\mathrm{S}^2)^{\times l N_{i,j}}.
\]
We have
\[
\gamma_1(P)|_C \cong \bigoplus_{\lambda \in S_1} \oplus_{m=1}^{lm(\lambda)} \lambda^*(\xi^{\otimes l}),
\]
\[
\gamma_2(P)|_C \cong \theta_lr,
\]
and
\[
\gamma_1(\Theta_2)|_C \cong \theta_{2 \mathrm{mult}(\gamma_1)},
\]
where $r \leq \mathrm{mult}(\gamma_2)$.  \cite[Lemma 2.1]{To2}
and (\ref{est1}) together imply that 
\[
\theta_{2 \mathrm{mult}(\gamma_1)} \precsim \left( \bigoplus_{\lambda \in S_1} \oplus_{m=1}^{lm(\lambda)} \lambda^*(\xi^{\otimes l}) \right) \oplus \theta_{lr}
\]
in the sense of Murray and von Neumann.  In other words, there is a $t \in \mathbb{N}$ and a complex
vector bundle $\omega$ over $C$ of fibre dimension $(l-2)\mathrm{mult}(\gamma_1) + lr$ such that
\[
\theta_{2 \mathrm{mult}(\gamma_1) + t} \oplus \omega \cong 
\left( \bigoplus_{\lambda \in S_1} \oplus_{m=1}^{lm(\lambda)} \lambda^*(\xi^{\otimes l}) \right) \oplus \theta_{lr+t}.
\]
Applying the total Chern class to this equation yields
\begin{eqnarray*}
c(\omega) & = & c\left(\left( \bigoplus_{\lambda \in S_1} \oplus_{m=1}^{lm(\lambda)} \lambda^*(\xi^{\otimes l}) \right) \oplus \theta_{lr+t} \right) \\
& = & \prod_{\lambda \in S_1} c(\lambda^*(\xi^{\otimes l}))^{l m(\lambda)} \\
& = & \prod_{\lambda \in S_1} [\lambda^*(c(\xi^{\otimes l}))]^{l m(\lambda)}. 
\end{eqnarray*}
Let us take the elements of $S_1$ to be numbered $\lambda_1,\ldots,\lambda_{N_{i,j}}$, so that
\[
c(\omega) = \prod_{k=1}^{N_{i,j}} (1 + e_{1,k} + \cdots + e_{l,k})^{l m(\lambda)}.
\]
It follows from our description of the ring structure of $H^*((\mathrm{S}^2)^{\times l})^{\otimes N_{i,j}}$
that we have $c_{lN_{i,j}}(\omega) \neq 0$, which in turn necessitates $\mathrm{rank}(\omega) \geq
lN_{i,j}$.  We conclude that
\begin{eqnarray*}
lN_{i,j} & \leq & (l-2)\mathrm{mult}(\gamma_1) + lr \\
& \leq  & (l-2)\mathrm{mult}(\gamma_1) + l \mathrm{mult}(\gamma_2) \\
& \leq & (l-2)M_{i,j} + 2(M_{i,j}-N_{i,j}).
\end{eqnarray*}
Dividing the last inequality above by $lM_{i,j}$ we get
\[
\frac{N_{i,j}}{M_{i,j}} \leq \frac{l-2}{l} + \frac{2}{l} \left(1- \frac{N_{i,j}}{M_{i,j}} \right).
\]
Using the assumption 
\[
\frac{N_{i,j}}{M_{i,j}} \geq \frac{6n-1}{6n} = \frac{2l-1}{2l}
\]
we have
\[
\frac{2l-1}{2l} \leq \frac{N_{i,j}}{M_{i,j}} \leq \frac{l-2}{l} + \frac{2}{l^2} < \frac{l-1}{l},
\]
a contradiction.
\end{proof}

\section{Strict comparison implies bounded dimension growth}\label{compbdg}

The next lemma says that if a simple $\vi$ algebra has strict comparison of positive elements, 
then it not only has slow dimension growth, but even
has slow dimension growth {\it as a $\vi$ algebra}.

\begin{lms}\label{sdg}
Let $A$ be a simple $\vi$ algebra;  
suppose that $A$ admits a standard decomposition
$(A_i,\phi_i)$ with seed space a CW-complex $X$. 
If $A$ has  strict comparison of positive elements,  
then $\mathrm{dim}(X)=0$ (in which case $A$ is AF), or for every $i \in \mathbb{N}$,
\[
\frac{N_{i,j}}{M_{i,j}} \stackrel{j \to \infty}{\longrightarrow} 0.
\]
If $A$ does not have strict comparison of positive elements, then
\begin{equation}
\label{N-M-limit-1}
\lim_{i \to \infty} \lim_{j \to \infty} \frac{N_{i,j}}{M_{i,j}} = 1.
\end{equation}
\end{lms}

\begin{proof}
Suppose, for a contradiction, that $A$ satisfies the hypotheses of
the lemma, $\mathrm{dim}(X) \geq 1$, and there is some $i_0 \in \mathbb{N}$
and $\delta > 0$ such that
\[
\frac{N_{i_0,j}}{M_{i_0,j}} \geq \delta, \ \forall j > i_0.
\]
We must show that $A$ does not have strict comparison of positive elements and 
that \eqref{N-M-limit-1} holds.

$A$ is simple with $\mathrm{dim}(X)=1$, so $M_{i,j} \to \infty$
as $j \to \infty$ for any fixed $i \in \mathbb{N}$.  This forces
$N_{i,j} \to \infty$, too, so that $\mathrm{dim}(X_i) \to \infty$
as $i \to \infty$. 

For any $j > m > i_0$ we have
\[
\delta \leq \frac{N_{i_0,j}}{M_{i_0,j}} = \frac{N_{i_0,m}}{M_{i_0,m}} \cdot
\frac{N_{m,j}}{M_{m,j}}.
\]
The sequence $(\frac{N_{i_0,j}}{M_{i_0,j}})_{j>i_{0}}$ is  decreasing, so its limit 
exists and is larger than or equal to $\delta$. It follows that $N_{m,j}/M_{m,j}$ 
approaches 1 as $m,j \to \infty$, whence \eqref{N-M-limit-1} holds.  Now apply Lemma \ref{ob} 
to see that $A$ must fail to have strict comparison of positive elements.
\end{proof}

\begin{props}\label{comparisonimpliesbdd}
Let $A$ be a simple $\vi$ algebra with strict comparison of positive elements, and suppose that $A$ admits
a standard decomposition $(A_i,\phi_i)$ with seed space a finite-dimensional CW-complex $X$.  
Then, $A$ has bounded dimension growth.
\end{props}

\begin{proof}
$A$ satisfies the hypotheses of Lemma \ref{sdg}.  If $A$ is AF, then it has bounded
dimension growth, so we may assume that $\mathrm{dim}(X) \geq 1$.  The conclusion of
Lemma \ref{sdg} then implies that
\begin{equation}\label{vanish}
\frac{N_{i,j}}{M_{i,j}} \stackrel{j \to \infty}{\longrightarrow} 0
\end{equation}
for every $i \in \mathbb{N}$.
We will prove that $A$ has very slow dimension growth in the sense of Gong;
bounded dimension growth follows by the Reduction Theorem of \cite{G} or, alternatively,
the classification theorem of \cite{Li2}.

Let there be given a positive tolerance $\epsilon >0$.  For natural numbers $j > i$
let $\pi_1,\ldots,\pi_{N_{i,j}}$ be the co-ordinate projection maps from $X_j = X_i^{n_j/n_i}$ to $X_i$
appearing as eigenvalue maps of $\phi_{i,j}$, and  
let $l_{i,j}$ be the number of eigenvalue maps of $\phi_{i,j}$ which are point evaluations.  
Since $A$ is simple and $\mathrm{dim}(X) \geq 1$, $l_{i,j_0} > 0$ for some $j_0 > i$.  
Straightforward calculation then shows that there are at least 
\[
l_{i,j_0} \cdot M_{j_0,j} = M_{i,j} \cdot \frac{l_{i,j_0}}{M_{i,j_0}}
\]
point evaluations in the map $\phi_{i,j}$.  Combining this with (\ref{vanish}) yields
\[
L_j := \left\lfloor \frac{M_{i,j}}{N_{i,j}} \cdot \frac{l_{i,j_0}}{M_{i,j_0}} \right\rfloor \stackrel{j \to \infty}{\longrightarrow} \infty.
\]
In other words, if one wants to partition the eigenvalue maps of $\phi_{i,j}$ which are point evaluations into
$N_{i,j}$ roughly equally sized multisets, then these multisets become arbitrarily large
as $j \to \infty$.  (We say multisets as some point evaluations may well be repeated.)
Assume that we have specified such a partition, let $S_1,\ldots,S_{N_{i,j}}$ denote
the multisets in the partition, and assume that $\mathrm{dim}(X_i)^3/L_j < \epsilon$.  
Each $f \in S_l$, $1 \leq l \leq N_{i,j}$, factors through
any of the co-ordinate projections $\pi_1,\ldots,\pi_{N_{i,j}}$.  Factor $f \in S_l$ as $f = \tilde{f} \circ \pi_l$,
where $\tilde{f}:X_i \to X_i$ has range equal to a point;  let $R_l$ be the multiset of all maps from $X_i$ to itself
obtained in this manner.  Let $t(l)$ be the number of copies of $\pi_l$ appearing among the eigenvalue
maps of $\phi_{i,j}$.  Put 
\[
B_i^{(l)} = \mathrm{M}_{m_i(t(l)+|R_l|)} \otimes \mathrm{C}(X_i),
\]
and observe that
\[
\mathrm{rank}(\mathbf{1}_{B_i^{(l)}}) = m_i(t(l)+|R_l|) \geq |R_l| = |S_l| \geq L_j.
\]
Define a map $\psi^{(l)}:A_i \to B_i^{(l)}$ by
\[
\psi^{(l)}(a) = \left( \bigoplus_{m=1}^{t(l)} a \right) \oplus \left( \bigoplus_{\tilde{f} \in R_l} a \circ \tilde{f} \right).
\]
Put $B_i = B_i^{(1)} \oplus \cdots \oplus B_i^{(N_{i,j})}$, and let $\psi:A_i \to B_i$
be the direct sum of the $\psi^{(l)}$, $1 \leq l \leq N_{i,j}$.  For each $1 \leq l \leq N_{i,j}$,
let $P_{l} \in A_j$ be the projection which is the sum of the images of the unit of $A_i$ under all
of the copies of $\pi_l$ in $\phi_{i,j}$ and all of the point evaluations in $S_l$.  Let 
$\gamma_l:B_i^{(l)} \to P_l A_j P_l$ be induced by $\pi_l$, and let $\gamma:B_i \to A_i$
be the direct sum of the $\gamma_l$.  We now have the factorisation
\[
A_i \stackrel{\psi}{\longrightarrow} B_i \stackrel{\gamma}{\longrightarrow} A_j,
\]
and each direct summand $B_i^{(l)}$ has the property that
\[
\frac{\mathrm{dim}(X_i)^3}{\mathrm{rank}(\mathbf{1}_{B_i^{(l)}})} \leq 
\frac{\mathrm{dim}(X_i)^3}{L_j} < \epsilon.
\]
Both $i$ and $\epsilon$ were arbitrary, so $A$ has very slow dimension growth.  
\end{proof}

Proposition \ref{comparisonimpliesbdd} establishes the implications $\mathrm{(ii)} \Rightarrow 
\mathrm{(iv)}$ and $\mathrm{(ii)} \Rightarrow \mathrm{(v)}$ of Theorem \ref{main}.  (The first observation
of the proof is that the hypotheses guarantee that $A$ has slow dimension growth as a
$\vi$ algebra, and so {\it a fortiori} as an AH algebra.)

\section{Finite decomposition rank}

In the present section we prove the remaining implication of Theorem \ref{main}, 
namely that finite decomposition rank implies strict comparison of positive elements 
in $\vi$ algebras. The technical key step is Lemma \ref{drcomp} below; under the 
additional assumption of real rank zero, a related result was already observed in 
\cite[Proposition 3.7]{Wi4}. Our proof is inspired by that argument.

\begin{lms}\label{drcomp}
Let $A$ be a simple, separable and unital $C^{*}$-algebra with $\dr A = n < \infty$. If $a,d^{(0)},\ldots,d^{(n)} \in A_{+}$ satisfy 
\[
d_{\tau}(a) < d_{\tau}(d^{(i)})
\]
for $i=0,\ldots,n$ and every $\tau \in T(A)$, then 
\[
\langle a \rangle \le \langle d^{(0)} \rangle + \ldots + \langle d^{(n)} \rangle.
\]
\end{lms}

\begin{proof}
It will be convenient to set up some notation: Given $0\le\alpha < \beta \le 1$, define  functions on the real line by
\[
g_{\beta}(t):= \left\{ 
\begin{array}{ll}
0 & t < \beta \\
1 & t \ge \beta 
\end{array}
\right., \; 
g_{\alpha,\beta}(t):= \left\{ 
\begin{array}{ll}
0 & t \le \alpha \\
1 & t \ge \beta \\
(t - \alpha)/(\beta - \alpha) & \mbox{ else}
\end{array}
\right.
\]
and
\[
f_{\alpha,\beta}(t):= \left\{ 
\begin{array}{ll}
0 & t \le \alpha \\
t & t \ge \beta \\
\beta (t - \alpha)/(\beta - \alpha) & \mbox{ else}.
\end{array}
\right. 
\]

Before turning to the proof, observe first that our hypotheses imply that $a$ is not invertible in $A$: indeed, if $a$ was invertible, we had $a^{1/n} \to \be_{A}$ as $n \to \infty$, so 
\[
1 = \lim_{i \to \infty} \tau(a^{1/n}) =  d_{\tau}(a) < d_{\tau}(d^{(i)}) \le 1,
\]
a contradiction, whence $0 \in \sigma(a)$. 

By passing to $M_{n+1}(A)$ (which again has decomposition rank $n$ by \cite{KW}, Corollary 3.9),  replacing each $d^{(i)} \in A \cong e_{00}M_{n+1}(A)e_{00}$ by a Cuntz equivalent element in the corner $e_{ii}M_{n+1}(A)e_{ii}$ for $i=1, \ldots,n$, and observing that $d_{\tau}(a) < d_{\tau}(d^{(i)})$ also holds for every $\tau \in T(M_{n+1}(A))$, we may as well  assume that the $d^{(i)}$ themselves are already pairwise orthogonal. 

For the actual proof of the lemma, we distinguish two cases. Suppose first that $0 \in \sigma(a)$ is an isolated point. Then, there is $\theta >0$ such that 
\[
p:= g_{\theta}(a) \in A
\]
is a projection satisfying $\langle p \rangle = \langle a \rangle$, whence $d_{\tau}(a) = d_{\tau}(p) = \tau(p)$ for all $\tau \in T(A)$. Furthermore, for any $\tau \in T(A)$ and $i=0, \ldots,n$ we have 
\begin{equation}
\label{dtauofdlimit}
d_{\tau}(d^{(i)}) = \lim_{\delta \searrow 0} \tau(g_{\delta/2,\delta}(d^{(i)})),
\end{equation}
so there are $\delta_{\tau}>0$ and $\eta_{\tau}>0$ such that 
\[
\tau(p) = d_{\tau}(a) < d_{\tau}(d^{(i)})  - \eta_{\tau} < \tau(g_{\delta_{\tau}/2,\delta_{\tau}}(d^{(i)}))
\]
for all $i$. Since the elements of $A$ are continuous when regarded as functions on $T(A)$, each $\tau$ has an open neighborhood $U_{\tau} \subset T(A)$ such that
\[
\tau'(p) < \tau'(g_{\delta_{\tau}/2,\delta_{\tau}}(d^{(i)})) 
\]
for $i=0, \ldots,n$ and $\tau' \in U_{\tau}$. Now by compactness of $T(A)$ (and since, for any positive $h$, 
\begin{equation}
\label{g-delta-monotony}
g_{\delta'/2,\delta'}(h) \le g_{\delta/2,\delta}(h)
\end{equation}
if only $\delta \le \delta'$) it is straightforward to find $\delta_{1}>0$ such that
\[
\tau(p) < \tau(g_{\delta_{1}/2,\delta_{1}}(d^{(i)}))
\]
for all $i=0, \ldots,n$ and  $\tau \in T(A)$. Now by \cite{Wi4}, Proposition 3.7, we have 
\[
p \precsim g_{\delta_{1}/2,\delta_{1}}(d^{(0)}) + \ldots + g_{\delta_{1}/2,\delta_{1}}(d^{(n)}),
\]
whence 
\[
\langle a \rangle = \langle p \rangle \le \langle g_{\delta_{1}/2,\delta_{1}}(d^{(0)}) \rangle + \ldots + \langle g_{\delta_{1}/2,\delta_{1}}(d^{(n)}) \rangle \le \langle d^{(0)} \rangle + \ldots + \langle d^{(n)} \rangle.
\]

Next, suppose $0$ is a limit point of $\sigma(a)$. The proof in this case is similar to that of \cite{Wi4}, Proposition 3.7, but we have to deal with some extra technical difficulties.  Since $a = \lim_{\epsilon \searrow 0} f_{\epsilon/2,\epsilon}(a)$ and since the $d^{(i)}$ are pairwise orthogonal, it will be enough to show that 
\[
f_{\epsilon/2,\epsilon}(a) \precsim d^{(0)} + \ldots + d^{(n)}
\]
for all $\epsilon>0$. So, given some $\epsilon>0$, we set 
\[
b:= f_{\epsilon/2,\epsilon}(a) \mbox{ and } c:= (g_{0,\epsilon/4} - g_{\epsilon/4,\epsilon/2})(a),
\]
then
\[
c \perp b \mbox{ and } b+c \precsim a.
\]
Since $0$ is a limit point of $\sigma(a)$, we have $c \neq 0$, hence (each $\tau \in T(A)$ is faithful by simplicity of $A$, $c$ is continuous as a function on $T(A)$ and $T(A)$ is compact)
\[
\alpha:= \min \{\tau(c) \, |\, \tau \in T(A)\} >0.
\]
Using that $c \le \be_{A}$ and that $c \perp b$, we obtain for all $\tau \in T(A)$
\begin{eqnarray*}
d_{\tau}(b) + \alpha & \le & d_{\tau}(b) + \tau(c) \\
& \le & d_{\tau}(b) + d_{\tau}(c) \\
& = & d_{\tau}(b+c) \\
&\le & d_{\tau}(a)
\end{eqnarray*}
and, by hypothesis,
\begin{equation}
\label{dtaubalpha}
d_{\tau}(b) < d_{\tau}(d^{(i)}) - \alpha
\end{equation}
for $i=0, \ldots,n$ and $\tau \in T(A)$.

Again, to show that $b \precsim d^{(0)}+ \ldots + d^{(n)}$ it will suffice to prove that
\[
f_{\eta,2\eta}(b) \precsim d^{(0)} + \ldots + d^{(n)}
\]
for any given $\eta>0$. To this end, we set
\begin{equation}
\label{barbdef}
\bar{b}:= g_{\eta/2,\eta}(b)
\end{equation}
and choose $0< \delta_{2}<\alpha/4$ such that
\begin{equation}
\label{taubarbalpha}
\tau(\bar{b})< \tau(\bar{d}^{(i)}) - \frac{3 \alpha}{4}
\end{equation}
for $i=0, \ldots, n$ and $\tau \in T(A)$, where
\[
\bar{d}^{(i)}:= g_{\delta_{2}/2,\delta_{2}}(d^{(i)}).
\]
The number $\delta_{2}$ is obtained in a similar way as $\delta_{1}$ in the first part of the proof, using compactness of $T(A)$: From \eqref{dtauofdlimit} and \eqref{dtaubalpha} we see that for each $\tau \in T(A)$ there is $\delta_{\tau}>0$ such that for $i=0, \ldots,n$
\[
\tau(\bar{b}) \le d_{\tau}(b) \stackrel{\eqref{dtaubalpha}}{<} d_{\tau}(d^{(i)})-\alpha \stackrel{\eqref{dtauofdlimit}}{\le} \tau(g_{\delta_{\tau}/2,\delta_{\tau}}(d^{(i)}))-\frac{3\alpha}{4}.
\]
Each $\tau$ has an open neighborhood $U_{\tau}$ such that
\[
\tau'(\bar{b}) < \tau'(g_{\delta_{\tau}/2,\delta_{\tau}}(d^{(i)})) - \frac{3 \alpha}{4}
\]
for $i=0, \ldots, n$ and $\tau' \in U_{\tau}$. Similar as in the first part of the proof, compactness of $T(A)$ and \eqref{g-delta-monotony} now yield $\delta_{2}>0$ such that \eqref{taubarbalpha} holds.

Since $\dr A=n$, by \cite{KW}, Proposition 5.1, there is a system $(F_{k},\psi_{k},\varphi_{k})_{k\in \mathbb{N}}$ of c.p.\ approximations for $A$ such that the $\varphi_{k}$ are $n$-decomposable and the $\psi_{k}$ are approximately multiplicative. In other words, for each $k \in \mathbb{N}$ there are finite-dimensional $C^{*}$-algebras $F_{k}$ and c.p.c.\ maps 
\[
A \stackrel{\psi_{k}}{\longrightarrow} F_{k} \stackrel{\varphi_{k}}{\longrightarrow} A
\]
such that 
\begin{enumerate}
\item $\varphi_{k} \psi_{k}(a) \to a$ for each $a \in A$ as $k \to \infty$ \label{pointwise-approximation}
\item $F_{k}$ admits a decomposition $F_{k}= \bigoplus_{i=0}^{n} F_{k}^{(i)}$ such that 
\[
\varphi_{k}^{(i)}:= \varphi_{k}|_{F_{k}^{(i)}}
\]
preserves orthogonality (i.e., has  order zero in the sense of  \cite{Wi1}, Definition 2.1b)) for each $i=0, \ldots,n$ and $k \in \mathbb{N}$
\item $\|\psi_{k}(aa') - \psi_{k}(a) \psi_{k}(a')\| \to 0$ for any $a,a' \in A$ as $k \to \infty$. \label{appr-mult-approximation}
\end{enumerate} 
We set 
\[
\psi_{k}^{(i)}(\, . \,):= \be_{F_{k}^{(i)}} \psi_{k}(\, . \, )
\]
for each $i$ and $k$ and note that the $\psi_{k}^{(i)}$ are also approximately multiplicative for each $i$ since the $\be_{F_{k}^{(i)}}$ are central projections in $F_{k}$. As in \cite{KW}, Remark 5.2(ii), we may (and will) assume that the $\psi_{k}$ are unital.  

Recall from \cite{Wi4}, 1.2, that each of the order zero maps $\varphi_{k}^{(i)}$ has a supporting $*$-homomorphism 
\[
\sigma_{k}^{(i)}:F_{k}^{(i)} \to A'';
\]
this a $*$-homomorphism satisfying 
\[
\varphi_{k}^{(i)}(x) = \sigma_{k}^{(i)}(x) \varphi_{k}^{(i)}(\be_{F_{k}^{(i)}}) = \varphi_{k}^{(i)}(\be_{F_{k}^{(i)}}) \sigma_{k}^{(i)}(x) \in A
\]
for all $x \in F_{k}^{(i)}$.

We proceed to show that there is $K \in \mathbb{N}$ such that
\begin{equation}
\label{taugalpha}
\tau( g_{1-\alpha/4}(\psi_{k}^{(i)}(\bar{b}))) < \tau(g_{\alpha/4}(\psi_{k}^{(i)}(\bar{d}^{(i)})))
\end{equation}
for all $i=0, \ldots, n$, $\tau \in T(F_{k}^{(i)})$ and $k \ge K$. If this was not the case, there would be s strictly increasing sequence $(k_{l})_{l \in \mathbb{N}} \subset \mathbb{N}$ such that, for some fixed $i_{0} \in \{0, \ldots,n\}$, there are  $\tau_{l} \in T(F_{k_{l}}^{(i_{0})})$ satisfying
\begin{equation}
\label{taugalphacontradiction}
\tau_{l}( g_{1-\alpha/4}(\psi_{k_{l}}^{(i_{0})}(\bar{b}))) \ge \tau_{l}(g_{\alpha/4}(\psi_{k_{l}}^{(i_{0})}(\bar{d}^{(i_{0})})))
\end{equation}
for all  $l\in \mathbb{N}$. But then 
\begin{eqnarray}
\tau_{l}(\psi_{k_{l}}^{(i_{0})}(\bar{b})) & \ge & \tau_{l}(g_{1-\alpha/4}(\psi_{k_{l}}^{(i_{0})}(\bar{b}))) - \frac{\alpha}{4} \nonumber \\
& \stackrel{\eqref{taugalphacontradiction}}{\ge} & \tau_{l}(g_{\alpha/4}(\psi_{k_{l}}^{(i_{0})}(\bar{d}^{(i_{0})}))) - \frac{\alpha}{4} \nonumber \\
& \ge & \tau_{l}(\psi_{k_{l}}^{(i_{0})}(\bar{d}^{(i_{0})})) - 2 \cdot \frac{\alpha}{4} \label{taubarbbard} \end{eqnarray}
for all $l \in \mathbb{N}$. Now fix some free ultrafilter $\omega \in \beta \mathbb{N} \setminus \mathbb{N}$, then
\[
\bar{\tau}(\, . \, ):= \lim_{\omega} \tau_{l}\psi_{k_{l}}^{(i_{0})}(\, . \,)
\]
obviously is a well-defined positive functional on $A$. It is tracial, since the $\tau_{l}$ are traces and the $\psi_{k_{l}}^{(i_{0})}$ are approximately multiplicative. It is a state, since $\bar{\tau}(\be_{A})=1$ (the $\tau_{l}$ are states and the $\psi_{k_{l}}^{(i_{0})}$ are unital). We have now constructed a tracial state  $\bar{\tau}$ on $A$ satisfying 
\[
\bar{\tau}(\bar{b}) \ge \bar{\tau}(\bar{d}^{(i_{0})}) - \frac{\alpha}{2},
\]
a contradiction to \eqref{taubarbalpha}. Therefore, there is $K \in \mathbb{N}$ such that \eqref{taugalpha} holds for all $i=0, \ldots, n$, $\tau \in T(F_{k}^{(i)})$ and $k \ge K$.

As a consequence, for $i=0,\ldots,n$ and $k\ge K$ there exist partial isometries $v_{k}^{(i)} \in F_{k}^{(i)}$ such that
\begin{equation}
\label{v*v}
(v_{k}^{(i)})^{*} v_{k}^{(i)} = g_{1-\alpha/4}(\psi_{k}^{(i)}(\bar{b})) \;  (\ge g_{1-\alpha/4,1}(\psi_{k}^{(i)}(\bar{b})))
\end{equation}
and 
\[
v_{k}^{(i)}(v_{k}^{(i)})^{*} \le g_{\alpha/4}(\psi_{k}^{(i)}(\bar{d}^{(i)}))
\]
(the $g_{1-\alpha/4}(\psi_{k}^{(i)}(\bar{b}))$ and $g_{\alpha/4}(\psi_{k}^{(i)}(\bar{d}^{(i)}))$ are projections in  $F_{k}^{(i)}$ --- which in turn are finite-dimensional algebras, hence satisfy the comparison property).

Note that
\begin{eqnarray*}
v_{k}^{(i)}(v_{k}^{(i)})^{*} & \le & g_{\alpha/4}(\psi_{k}^{(i)}(\bar{d}^{(i)})) \\
& \le & g_{\alpha/4}(\psi_{k}(\bar{d}^{(i)})) \\
& \le & g_{0,\delta_{2}}(\psi_{k}(\bar{d}^{(i)})) \\
& \le & \be_{F_{k}}
\end{eqnarray*}
for $i=0, \ldots,n$ and  $k \ge K$; since the $v_{k}^{(i)}(v_{k}^{(i)})^{*}$ are projections, from this one easily concludes that 
\[
g_{0,\delta_{2}}(\psi_{k}( \bar{d}^{(i)})) v_{k}^{(i)}(v_{k}^{(i)})^{*} = v_{k}^{(i)}(v_{k}^{(i)})^{*}.
\]
Because the $\psi_{k}$ are approximately multiplicative, and using that the $\bar{d}^{(i)}$ are mutually orthogonal, we also have
\[
\|g_{0,\delta_{2}}(\psi_{k}(\sum_{j=0}^{n} \bar{d}^{(j)})) - \psi_{k}(\sum_{j=0}^{n} g_{0,\delta_{2}}(\bar{d}^{(j)}))\| \stackrel{k \to \infty}{\longrightarrow} 0
\]
and
\[
\|g_{0,\delta_{2}}(\psi_{k}(\bar{d}^{(i')})) v_{k}^{(i)}(v_{k}^{(i)})^{*} - \delta_{i',i} \cdot v_{k}^{(i)}(v_{k}^{(i)})^{*}\| \stackrel{k \to \infty}{\longrightarrow}
\]
for $i,i' \in \{0, \ldots,n\}$ ($\delta_{i,i'}$ denotes the Kronecker delta of $i$ and $i'$).  
Moreover, we have 
\[
\|\sum_{j=0}^{n} g_{0,\delta_{2}}(\bar{d}^{(j)}) \varphi_{k}(v_{k}^{(i)}(v_{k}^{(i)})^{*}) - \varphi_{k}(\psi_{k}(\sum_{j=0}^{n} g_{0,\delta_{2}}(\bar{d}^{(j)})) v_{k}^{(i)}(v_{k}^{(i)})^{*})\| \le 2 \mu_{k}^{\frac{1}{2}}
\]
by \cite{KW}, Lemma 3.5 (an easy consequence of Stinespring's Theorem), where
\[
\mu_{k}:= \max \{\|(\varphi_{k}\psi_{k} - \id)(\sum_{j=0}^{n} g_{0,\delta_{2}}(\bar{d}^{(j)}))\| , \, \|(\varphi_{k}\psi_{k} - \id)((\sum_{j=0}^{n} g_{0,\delta_{2}}(\bar{d}^{(j)}))^{2})\| \}.
\]
Observing that $\mu_{k} \to 0$ as $k \to \infty$ and combining all these facts we obtain
\[
\|\sum_{j=0}^{n} g_{0,\delta_{2}}(\bar{d}^{(j)}) \varphi_{k}(v_{k}^{(i)} (v_{k}^{(i)})^{*}) - \varphi_{k}(v_{k}^{(i)} (v_{k}^{(i)})^{*}) \| \stackrel{k \to \infty}{\longrightarrow} 0
\]
and
\[
\|g_{0,\delta_{2}}( \bar{d}^{(i')}) \varphi_{k}(v_{k}^{(i)} (v_{k}^{(i)})^{*}) - \delta_{i,i'} \cdot \varphi_{k}(v_{k}^{(i)} (v_{k}^{(i)})^{*}) \| \stackrel{k \to \infty}{\longrightarrow} 0
\]
for $i,i'=0, \ldots, n$. Using that 
\[
\varphi_{k}(v_{k}^{(i)}(v_{k}^{(i)})^{*}) = \varphi_{k}^{(i)}(v_{k}^{(i)}(v_{k}^{(i)})^{*}) = \varphi_{k}^{(i)}(\be_{F_{k}^{(i)}}) \sigma_{k}^{(i)}(v_{k}^{(i)})\sigma_{k}^{(i)}(v_{k}^{(i)})^{*}
\]
it follows easily that 
\begin{eqnarray*}
\|\sigma_{k}^{(i')}(v_{k}^{(i')})^{*}\varphi_{k}^{(i')}(\be_{F_{k}^{(i')}})^{\frac{1}{2}} \sum_{j=0}^{n} g_{0,\delta_{2}}(\bar{d}^{(j)}) \varphi_{k}^{(i)}(\be_{F_{k}^{(i)}})^{\frac{1}{2}} \sigma_{k}^{(i)}(v_{k}^{(i)}) \\
- \delta_{i',i} \cdot \varphi_{k}^{(i)}((v_{k}^{(i)})^{*}v_{k}^{(i)}) \| \stackrel{k \to \infty}{\longrightarrow} 0 
\end{eqnarray*}
for $i,i'=0, \ldots, n$, whence 
\[
\|s_{k}^{*}\sum_{j=0}^{n} g_{0,\delta_{2}/2}(\bar{d}^{(j)}) s_{k} - f_{\eta,2\eta}(b)^{\frac{1}{2}}(\sum_{i=0}^{n} \varphi_{k}^{(i)}((v_{k}^{(i)})^{*} v_{k}^{(i)})) f_{\eta,2\eta}(b)^{\frac{1}{2}}\| \stackrel{k \to \infty}{\longrightarrow} 0,
\]
where
\[
s_{k}:= \sum_{i=0}^{n} \varphi_{k}^{(i)}(\be_{F_{k}^{(i)}})^{\frac{1}{2}} \sigma_{k}^{(i)}(v_{k}^{(i)}) f_{\eta,2 \eta}(b)^{\frac{1}{2}}
\]
for $k \ge K$. We now have 
\begin{eqnarray*}
\lefteqn{\lim_{k \to \infty} \|f_{\eta,2\eta}(b) - s_{k}^{*}\sum_{j=0}^{n} g_{0,\delta_{2}}(\bar{d}^{(j)})s_{k} \| }\\
& = & \lim_{k \to \infty} \|f_{\eta,2\eta}(b) - f_{\eta,2\eta}(b)^{\frac{1}{2}} (\sum_{i=0}^{n} \varphi_{k}^{(i)}((v_{k}^{(i)})^{*} v_{k}^{(i)})) f_{\eta,2\eta}(b)^{\frac{1}{2}}  \| \\
& \stackrel{\eqref{v*v}}{\le} & \lim_{k \to \infty} \|f_{\eta,2\eta}(b) - f_{\eta,2\eta}(b)^{\frac{1}{2}} (\sum_{i=0}^{n} \varphi_{k}^{(i)}(g_{1-\alpha/4,1}(\psi_{k}^{(i)}(\bar{b})))) f_{\eta,2\eta}(b)^{\frac{1}{2}}  \| \\
& \stackrel{\mbox{\footnotesize \ref{appr-mult-approximation}}}{ =} & \lim_{k \to \infty} \|f_{\eta,2\eta}(b) - f_{\eta,2\eta}(b)^{\frac{1}{2}} (\sum_{i=0}^{n} \varphi_{k}^{(i)}\psi_{k}^{(i)}(g_{1-\alpha/4,1}(\bar{b}))) f_{\eta,2\eta}(b)^{\frac{1}{2}}  \| \\
& \stackrel{\mbox{\footnotesize \ref{pointwise-approximation}}}{=} & \lim_{k \to \infty} \|f_{\eta,2\eta}(b) - f_{\eta,2\eta}(b)^{\frac{1}{2}} g_{1-\alpha/4,1}(\bar{b}) f_{\eta,2\eta}(b)^{\frac{1}{2}}  \| \\
& \stackrel{\eqref{barbdef}}{=} & 0.
\end{eqnarray*}
This shows that
\[
f_{\eta,2\eta}(b) \precsim \sum_{j=0}^{n} g_{0,\delta_{2}}(\bar{d}^{(j)}).
\]
Since $\eta$ was arbitrary, and because 
\[
\sum_{j=0}^{n} g_{0,\delta_{2}}(\bar{d}^{(j)}) \precsim \sum_{j=0}^{n}\bar{d}^{(j)},
\]
it follows that 
\[
b \precsim \sum_{j=0}^{n} d^{(j)},
\]
as desired.
\end{proof}

\begin{cors}\label{finitedrimpliescomparison}
Let $A$ be as in the hypotheses of Theorem \ref{main}.  If $A$ has finite decomposition
rank, then $A$ has strict comparison of positive elements.
\end{cors}

\begin{proof}
We prove the contrapositive.  Suppose that $A$ does not have strict comparison of positive elements and fix a standard decomposition as in \ref{main}.  Then, condition \eqref{N-M-limit-1} of Lemma \ref{sdg} holds.  It follows that $A$  satisfies the hypotheses of
Lemma \ref{ob}, and so strict comparison fails in the manner prescribed the in the conclusion of that Lemma.
In light of Lemma \ref{drcomp}, this failure excludes the possibility that $A$ has finite decomposition rank.
\end{proof}

The preceding corollary establishes the implication (iii)$\Rightarrow$(ii) of Theorem \ref{main}.

\section{Real rank zero}

In this section we prove that an algebra of real rank zero which also satisfies 
the hypotheses of Theorem \ref{main} must then satisfy conditions (i)-(vi) of the same theorem, thus completing the proof of our main result.  
The result is a special case of a theorem of Toan Ho and the first named author which will appear in
Toan Ho's Ph.D. thesis.  As no preprint of this result was available at the time of writing,
we give a proof here which applies only to $\vi$ algebras.

Let $X$ be a compact connected Hausdorff space and $a$ a self-adjoint element of $\mathrm{M}_n(\mathrm{C}(X))$.
For each $x \in X$, form an $n$-tuple consisting of the eigenvalues of $a$ listed in decreasing order.  
For each $m \in \{1,\ldots,n\}$ let $\lambda_m:X \to \mathbb{R}$ be the function whose value at $x$
is the $m^{\mathrm{th}}$ entry of the eigenvalue $n$-tuple for $x$.  
The {\it variation of the normalised trace} of $a$ (v. \cite{BBEK}), denoted $TV(a)$, is defined
as
\[
\mathrm{sup} \left\{ \left| \frac{1}{n} \sum_{m=1}^n (\lambda_m(x)-\lambda_m(y)) \right|  \ : \ x,y \in X \right\}.
\]
Suppose that $A = \lim_{i \to \infty}(\mathrm{M}_{m_i}(\mathrm{C}(X_i),\phi_i)$ is of real rank zero, and let $a$ be a self-adjoint
element of some $\mathrm{M}_{m_i}(\mathrm{C}(X_i))$.  Then, by Theorem 1.3 of \cite{BBEK}, the variation of the
normalised trace tends to zero as $j \to \infty$ for each direct summand of $\phi_{i,j}(a)$ corresponding to 
a connected component of  $X_{j}$.

\begin{props}
Let $A=\lim_{i \to \infty}(A_i,\phi_i)$ be a simple $\vi$ algebra with seed space a finite-dimensional CW-complex.
If $A$ has real rank zero, then $A$ has bounded dimension growth.
\end{props}

\begin{proof}
If $A$ is AF, then there is nothing to prove.  If $A$ is not AF, then all but finitely many of the $\phi_i$s
contain at least one co-ordinate projection as an eigenvalue map, and each $X_i$ has dimension strictly greater
than zero.  It will be enough to prove that for each $i \in \mathbb{N}$,
\[
\frac{N_{i,j}}{M_{i,j}} \stackrel{j \to \infty}{\longrightarrow} 0.
\]
The proof of Proposition \ref{comparisonimpliesbdd} then shows that $A$ has bounded dimension growth.

Let $\epsilon > 0$ be given, and
suppose for a contradiction that for some $i \in \mathbb{N}$ and $c > 0$ we have
\[
\frac{N_{i,j}}{M_{i,j}} \stackrel{j \to \infty}{\longrightarrow} c.
\]
By increasing $i$ if necessary (and following the lines of the proof of Lemma \ref{sdg}) we may assume that $c > 7/8$.
Choose a continuous
function $f:X_i \to [0,1]$ such that for some points $x_0,x_1$ in the same connected component
of $X_i$ we have
$f(x_0)=0$ and $f(x_1)=1$; put $a := f \cdot \mathbf{1}_{A_i}$.  

For any
$j > i$ we have
\[
\phi_{i,j}(a)(x) = \mathrm{diag}(a(\gamma_1(x)),\ldots,a(\gamma_{M_{i,j}}(x))), \ \forall x \in X_j,
\]
where the $\gamma_l$s are the eigenvalue maps of $\phi_{i,j}$.  Let 
$\pi_1,\ldots,\pi_{N_{i,j}}:X_j \to X_i$ be the distinct co-ordinate projections
appearing among the $\gamma_l$s.  Since $TV(\phi_{i,j}(a))$ is unaffected by unitary
conjugation in $A_j$, we may assume that
\[
\phi_{i,j}(a)(x) = \mathrm{diag}(a(\pi_1(x)),\ldots,a(\pi_{N_{i,j}}(x)),\ldots,
a(\gamma_{M_{i,j}}(x))), \ \forall x \in X_j.
\]  
Fix a point $y_0 \in X_j$ which when viewed as an element of a Cartesian power of $X_i$ has
the value $x_0$ in each co-ordinate;  define
$y_1$ similarly with respect to $x_1$, and notice that $y_0$ and $y_1$ are in the same connected 
component of $X_j$.  Then, the eigenvalue list of $\phi_{i,j}(a)(y_0)$ contains
at least $m_{i} N_{i,j}$ 0s, while the list for $\phi_{i,j}(a)(y_1)$ contains at least $m_{i} N_{i,j}$ 1s.
By the pigeonhole principle, at least $m_i[M_{i,j} - 2(M_{i,j}-N_{i,j})]$ of the eigenfuntions
$\lambda_m$ corresponding to $\phi_{i,j}(a)$ have the value 0 at $y_0$ and 1 at $y_1$, while the
remaining $2m_i(M_{i,j}-N_{i,j})$ eigenfunctions satisfy $\lambda_m(y_1)-\lambda_m(y_0) \geq -1$.  Now,
\begin{eqnarray*}
TV(\phi_{i,j}(a))
& \geq & \left| \frac{1}{m_i M_{i,j}} \sum_{m=1}^{m_j} (\lambda_m(y_1)-\lambda_m(y_0)) \right| \\
& \geq & \frac{M_{i,j} - 4(M_{i,j}-N_{i,j})}{M_{i,j}} \\
& = & \frac{4 N_{i,j}}{M_{i,j}} - 3 \\
& > &  \frac{1}{2}
\end{eqnarray*}  
since $N_{i,j}/M_{i,j}> c > 7/8$.   
This contradicts our real rank zero assumption for $A$, completing the proof.
\end{proof}

%%%%%%%%%%%%%%%%%%%%%%%%%%%%%%%%%%%%%%%%%%%%%%%

\section{Non-$\mathcal{Z}$-stable $\vi$ algebras}

We now give examples of non-isomorphic $\vi$-algebras which 
cannot be distinguished using topological $\mathrm{K}$-theory and traces.  These
are not the first such --- examples are already given in \cite{To2} --- but the results of the
present paper allow us to construct a large class of examples with relatively little further
effort.  They will also demonstrate the variety of tracial state spaces which can occur in a simple nuclear
C$^*$-algebra of infinite decomposition rank.  

Any subgroup $G$ of $\mathbb{Q}$ corresponds to a list of prime powers $P_G = \{p_1^{n_1},p_2^{n_2},\ldots\}$, $n_i \in \mathbb{Z}^+ \cup \{\infty\}$,
in the following sense: the elements of $G$ are those rationals which, when in loweset terms, have denominators
of the form $p_1^{r_1} p_2^{r_2} \cdots$, where $r_i < n_i + 1$ for all $i$ and $r_i=0$ for all but finitely many $i$.
If $n_i = \infty$ for some $i$ then we will say that $G$ is of {\it infinite type}.  Let $p$ be a prime.
If $p^{\infty} \in P_G$ and $H$ is the subgroup of $\mathbb{Q}$ with $P_H = \{p^{\infty}\}$, then
$H \otimes G \cong G$.   

Let $X$ be a contractible and finite-dimensional CW-complex.  Construct a $\vi$ algebra $A_X = \lim_{i \to \infty}(A_{i},\phi_i)$
satisfying:
\begin{enumerate}
\item[(i)] The ratio $N_{1,j}/M_{1,j}$ does not vanish;
\item[(ii)] $A_X$ is simple by virtue of a judicious inclusion of point evaluations as eigenvalue maps of the $\phi_i$;
\item[(iii)] The $\mathrm{K}_0$-group of $A_X$ (necessarily a subgroup of $\mathbb{Q}$ by the contractibility of $X_1 = X$
and each $X_i$) is of infinite type.
\end{enumerate}

\noindent
Inspection of Villadsen's  construction in \cite{V1} shows that for a fixed $X$, one can arrange for $\mathrm{K}_0(A_X)$ to be
an arbitrary infinite type subgroup of $\mathbb{Q}$.  There are uncountably many such subgroups, and hence,
for a fixed $X$, uncountably many non-isomorphic algebras $A_X$ satisfying $\mathrm{(i)}-\mathrm{(iii)}$.
Condition $\mathrm{(i)}$ ensures that $A_X$ does not have strict comparison of positive elements (use Lemmas
\ref{ob} and \ref{sdg}).

Fix an algebra $A_X$ as above.  Let $p$ be a prime such that $p^{\infty} \in P_{\mathrm{K}_0(A_X)}$, and let
$\mathfrak{U}$ be a UHF algebra with $P_{\mathrm{K}_0(\mathfrak{U})} = \{p^{\infty}\}$.  We claim that the tensor product
$A_X \otimes \mathfrak{U}$ has the same topological $\mathrm{K}$-theory and tracial state space as $A$.
At the level of $\mathrm{K}$-theory this statement follows from the K\"{u}nneth Theorem, the triviality of the
$\mathrm{K}_1$-groups of both $\mathfrak{U}$ and $A_X$ (in the case of $A_X$ this is due to the contractibility
of $X_i$), and the isomorphism
\[
\mathrm{K}_0(A_X) \otimes \mathrm{K}_0(\mathfrak{U}) \cong \mathrm{K}_0(A_X).
\]
At the level of tracial state spaces the statement is due to the fact that $\mathfrak{U}$ admits a
unique tracial state.  There is only one possible pairing of traces with $\mathrm{K}_0$
in each of $A_X$ and $A_X \otimes \mathfrak{U}$, as their $\mathrm{K}_0$-groups are subgroups of the
rationals.  As noted above, $A_X$ does not have strict comparison of positive elements, but 
$A_X \otimes \mathfrak{U}$ does by virtue of \cite[Lemma 5.1]{R6}.  Thus, $A_X$ and $A_X \otimes
\mathfrak{U}$ are not isomorphic, and by varying $X$ and $\mathrm{K}_0(A_X)$ independently 
we obtain a large class of examples of the desired variety. 

Straightforward but laborious calculation shows that the tracial state space of $A_X$ as above
is a Bauer simplex with extreme boundary homeomorphic to $X^{\times \infty}$.  (The details of this
calculation are more or less contained in the proof of \cite[Theorem 4.1]{To7} --- we will not reproduce
them here.)  $A_X$ also has infinite decomposition rank by Theorem \ref{main}, and this does not depend
on $X$ being contractible.  Thus, a large variety of structure can occur in the tracial state
space of a simple nuclear C$^*$-algebra with infinite decomposition rank.

Finally, we remark that infinite decomposition rank can also occur in the case of a simple AH algebra with \emph{unique} tracial state, as observed in \cite[Example 6.6(i)]{Wi4}, using the examples of \cite{V2}.


\begin{thebibliography}{999}


%\bibitem{Bl1} Blackadar, B.: $K$-Theory for Operator Algebras, 2nd edition, MSRI Publications {\bf 5}, Cambridge University Press 1998

%\bibitem{Bl2} Blackadar, B.: {\it Matricial and ultramatricial topology}, Operator algebras, mathematical physics, and low-dimensional
%topology (Istanbul, 1991), Res. Notes Math. {\bf 5}, 11-38

\bibitem{BBEK} Blackadar, B., Bratteli, O., Elliott, G. A., and Kumjian, A.: {Reduction of real rank in inductive limits of $C^*$-algebras}, 
Math.\ Ann.\ {\bf 292} (1992), 111-126

%\bibitem{BDR} Blackadar, B., Dadarlat, M., and R{\o}rdam, M.: {\it The
%real rank of inductive limit $C^{*}$-algebras}, Math. Scand. {\bf 69} 
%(1991), 211-216

\bibitem{BH} Blackadar, B., and Handelman, D.: {\it Dimension functions and traces on $C^*$-algebras}, J.\ Funct.\ Anal.\ {\bf 45} (1982), 297-340

\bibitem{BKR} Blackadar, B., Kumjian, A., and R{\o}rdam, M.: {\it Approximately central matrix units and the structure of noncommutative tori}, $\mathrm{K}$-theory {\bf 6} (1992), 267-284


\bibitem{Bra:AF} Bratteli, O.: {\it Inductive limits of finite dimensional C*-algebras}, Trans.\ Amer.\ Math.\ Soc.\ {\bf 171} (1972), 195-234

%\bibitem{BP} Brown, L.\ G.\ and Pedersen, G.\ K.: {\it $C^*$-algebras of real rank zero}, J.\ Funct.\ Anal.\ {\bf 99} (1991), 131-149

%\bibitem{B} Brown, N.: {\it Invariant means and finite representation theory of $C^{*}$-algebras}, preprint math.OA/0304009 (2003)

%\bibitem{CE} Choi, M. D., and Effros, E. G.: {\it The completely positive lifting problem for $C^*$-algebras}, Ann. of Math. (2) {\bf 104} (1976), 585-609

%\bibitem{C} Cuntz, J.: {\it Dimension functions on simple $C^*$-algebras}, Math. Ann. {\bf 233} (1978), 145-153

%\bibitem{DNNP} Dadarlat, M., Nagy, G., Nemethi, A., and Pasnicu, C.: {\it Reduction of topological stable rank in inductive limits of $C^*$-algebras},
%Pac. J. Math. {\bf 153} (1992), 267-276

%\bibitem{D} Dadarlat, M.: {\it Reduction to dimension three of of local spectra of real rank zero $C^{*}$-algebras}, J.\ Reine Angew.\ Math.\ {\bf 460} (1995), 189-212

%\bibitem{DG} Dadarlat, M.\ and Gong., G.: {\it A classification result for
%approximately homogeneous $C^*$-algebras of real rank zero}, Geom.\ Funct.\
%Anal.\ {\bf 7} (1997), 646-711

%\bibitem{DL} Dadarlat, M.\ and Loring, T.: {\it Classifying $C^*$-algebras
%via ordered, mod-p $\mathrm{K}$-theory}, Math.\ Ann.\ {\bf 305} (1996), 601-616

%\bibitem{Da} Davidson ???

%\bibitem{Dx} Dixmier, J.: Les $C^*$-alg\`ebres et leurs repr\'esentations, Gauthiers-Villars, Paris 1964

%\bibitem{ER} Effros, E.\ G.\ and Rosenberg, J.: {\it $C^{*}$-algebras with approximately inner flip}, Pac.\ J.\ Math.\ {\bf 77} (1978), 417-443

\bibitem{Ell:AF} Elliott, G.\ A.: {\it On the classification of inductive 
limits of sequences of semi-simple finite-dimensional algebras}, J.\ Algebra 
{\bf 38} (1976), 29-44

\bibitem{El2} Elliott, G.\ A.: {\it The classification problem for amenable
$C^*$-algebras}, Proc. ICM '94, Zurich, Switzerland, Birkhauser Verlag, Basel,
Switzerland, 922-932
 
\bibitem{El4} Elliott, G.\ A.: {\it An invariant for simple C*-algebras}, Canadian Mathematical Society. 
1945-1995, vol. 3, Canadian Math.\ Soc., Ottawa, ON, 1996, pp.\ 61-90
 
\bibitem{Ell:AT} Elliott, G.\ A.: {\it On the classification of $C^*$-algebras
of real rank zero}, J.\ Reine Angew.\ Math.\ {\bf 443} (1993), 179-219

%\bibitem{EE} Elliott, G.\ A.\ and Evans, D.\ E.: {\it The structure of irrational rotation 
%$C^*$-algebras }, Ann.\ Math.\ {\bf 138} (1993), 477-501

%\bibitem{EG} Elliott, G. A., and Gong, G.:{On the classification of $C^*$-algebras of real rank zero. II.} Ann. of Math. (2) {\bf 144} (1996), 497-610

\bibitem{EGL2} Elliott, G.\ A., Gong, G.\ and Li, L.: {\it Approximate divisibility of simple inductive limit $C^{*}$-algebras}, Contemp.\ Math.\ {\bf 228} (1998), 87-97


\bibitem{EGL} Elliott, G.\ A., Gong, G.\ and Li, L.: {\it On the classification
of simple inductive limit $C^*$-algebras, II:  The isomorphism theorem}, preprint

%\bibitem{EGL3} Elliott, G. A., Gong, G., and Li, L.: {\it Injectivity of the connecting maps in AH inductive limit systems}, Canad. Math. Bull. {\bf 48} (2005), 50-68

%\bibitem{ER} Elliott, G. A. and R{\o}rdam, M.: {\it Perturbation of Hausdorff moment sequences, and an application to the theory of $C^*$-algebras of real rank zero},
%preprint math.OA/0502239 (2005)

%\bibitem{EV} Elliott, G. A. and Villadsen, J.: {\it Perforated ordered $\mathrm{K}_0$-groups}, Canad. J. Math. 
%{\bf 52} (2000), 1164-1191

%\bibitem{Fe} Fell, J.\ M.\ G.: {\it The structure of algebras of operator fields}, Acta Math.\ {\bf 106} (1961), 233-280

\bibitem{Gli:UHF} Glimm, J.: {\it On a certain class of operator algebras}, Trans.\ Amer.\ Math.\ Soc.\ {\bf 95} (1960), 318-340

%\bibitem{Go1} Goodearl, K. R.: {\it Riesz decomposition in inductive limit $C^*$-algebras}, Rocky Mt. J. 
%Math. {\bf 24} 1994, 1405-1430

%\bibitem{Go2} Goodearl, K.\ R.: Partially Ordered Groups with Interpolation, Mathematical Surveys and Monographs {\bf 20}, AMS, Providence, RI, 1986

\bibitem{G} Gong, G.: {\it On the classification of simple inductive limit $C^*$-algebras I.  The reduction theorem.}, Doc.\ Math.\ {\bf 7} (2002), 255-461

%\bibitem{GJS} Gong, G., Jiang, X.\ and  Su, H.: {\it Obstructions to $\mathcal{Z}$-stability for unital simple $C^{*}$-algebras}, 
%Canad.\ Math.\ Bull.\ {\bf 43} (2000), 418-426

%\bibitem{Ha} Haagerup, U.: {\it Every quasi-trace on an exact $C^*$-algebra is a trace}, preprint, 1991

%\bibitem{H} Husemoller, D.: Fibre Bundles, McGraw-Hill, New York 1966

%\bibitem{HW} Hurewicz, W.\ and Wallmann, H.: Dimension Theory, Princeton University Press, Princeton 1941

%\bibitem{J} Jiang, X.: {\it Nonstable $K$-theory for $\mathcal{Z}$-stable $C^{*}$-algebras}, preprint math.OA/9707228 (1997)

\bibitem{JS1} Jiang, X.\ and Su, H.: {\it On a simple unital projectionless 
$C^{*}$-algebra}, Amer.\ J.\ Math.\ {\bf 121} (1999), 359-413

\bibitem{Kir:pi} Kirchberg, E.: {\it The classification of Purely Infinite $C^*$-algebras using Kasparov's Theory}, Fields Inst. Comm., to appear

\bibitem{KW} Kirchberg, E.\ and Winter, W.: {\it Covering dimension and 
quasidiagonality}, Intern.\ J.\ Math.\ {\bf 15} (2004), 63-85

%\bibitem{Li1} Lin, H.: {\it Tracially AF $C^{*}$-algebras},  Trans.\ Amer.\ Math.\ Soc.\  {\bf 353}  (2001), 693-722

%\bibitem{Li3} Lin, H.: {\it The tracial topological rank of $C^*$-algebras}, Proc. London Math. Soc. (3) {\bf 83} (2001), 199-234

%\bibitem{Lin:TAF} Lin, H.:  {\it Classification of simple tracially AF 
%$C^{*}$-algebras},  Canad.\ J.\ Math.\  {\bf 53}  (2001),  161-194

\bibitem{Li2} Lin, H.: {\it Simple nuclear $C^*$-algebras of tracial
topological rank one}, arXiv preprint math.OA/0011079 (2004)

\bibitem{Lin:tr0} Lin, H.: {\it Classification of simple C*-algebras with tracial topological
rank zero}, Duke Math.\ J.\ {\bf 125} (2004), 91-119


%\bibitem{Li4} Lin, H.: {\it Simple AH algebras of real rank zero}, Proc. Amer. Math. Soc. {\bf 131} (2003), 3813-3819

%\bibitem{M} Mygind, J.: {\it Classification of certain simple $C^*$-algebras with torsion in $\mathrm{K}_1$},
%Canad.\ J.\ Math.\ {\bf 53} (2001), 1223-1308

%\bibitem{N} Nistor, V.: {\it Stable rank for a certain class of type I $C^*$-algebras}, J. Operator Theory {\bf 17} (1987), 365-373

%\bibitem{Pa} Paulsen, V.: Completely bounded maps and dilations, Pitman Research Notes in Mathematics {\bf 146} (1986)

%\bibitem{Pe} Pears, A.\ R.: Dimension Theory of General Spaces, Cambridge Univ. Press, Cambridge 1975

%\bibitem{P} Perera, F.:{\it Monoids arising from positive matrices over commutative $C^*$-algebras}, Math. Proc. R. Ir. Acad. {\bf 99A} (1999), 75-84

%\bibitem{RP} Perera, F. and R{\o}rdam, M.: {\it AF-embeddings into $C^{*}$-algebras with real rank zero}, J. Funct. Anal. {\bf 217} (2004), 142-170

\bibitem{Phi:pi} Phillips, N.\ C.: {\it A classification theorem for nuclear
purely infinte simple $C^*$-algebras}, Doc.\ Math.\ {\bf 5} (2000), 49-114

%\bibitem{P3} Phillips, N.\ C.: {\it Cancellation and stable rank for direct limits of recursive subhomogeneous algebras}, arXiv preprint math.OA/0101157 (2001)

%\bibitem{P4} Phillips, N.\ C.: {\it Recursive subhomogeneous algebras}, Trans.\ Amer.\ Math.\ Soc., to appear

%\bibitem{P3} Phillips, N.\ C.: {\it Recursive subhomogeneous algebras}, preprint math.OA/0101156 (2001)

%\bibitem{Ri1} Rieffel, M.\ A.: {\it Dimension and stable rank in the $K$-theory of $C^*$-algebras}, Proc.\ London Math.\ Soc.\ (3) {\bf 46} 
%(1983), 301-333

\bibitem{R1} R{\o}rdam, M.: Classification of Nuclear $C^*$-Algebras, 
Encyclopaedia of Mathematical Sciences {\bf 126}, Springer-Verlag, Berlin, 
Heidelberg 2002

\bibitem{R3} R{\o}rdam, M.: {\it A simple $C^*$-algebra with a finite and an
infinite projection}, Acta Math.\ {\bf 191} (2003), 109-142

\bibitem{R4} R{\o}rdam, M.: {\it The stable and the real rank of $\mathcal{Z}$-absorbing $C^{*}$-algebras}, Int.\ J.\ Math.\ {\bf 15} (2004), 1065-1084

%\bibitem{R5} R{\o}rdam, M.: {\it Stability of $C^*$-algebras is not a stable property}, Doc.\ Math.\  {\bf 2} (1997), 375-386

%\bibitem{Ror:gra} R\o rdam, M.: {\it Classification of Cuntz-Krieger algebras}, $\mathrm{K}$-theory {\bf 9} (1995), 31-58

\bibitem{R6} R\o rdam, M.: {\it On the structure of simple C$^*$-algebras tensored with a UHF-algebra, II}, J. Funct. Anal. {\bf 107}
(1992), 255-269

\bibitem{To1} Toms, A.\ S.: {\it On the independence of $\mathrm{K}$-theory
and stable rank for simple $C^*$-algebras}, J.\ Reine Angew.\ Math. {\bf 578} (2005), 185-199

\bibitem{To2} Toms, A.\ S.: {\it On the classification problem for nuclear $C^*$-algebras}, arXiv preprint math.OA/0509103 (2005);  
Ann.\ of Math.\ (2), to appear


%\bibitem{To3} Toms, A.\ S.: {\it Dimension growth for $C^*$-algebras}, arXiv preprint math.OA/0509159 (2005)

%\bibitem{To4} Toms, A.\ S.: {\it Cancellation does not imply stable rank one}, Bull. London Math. Soc., to appear

\bibitem{To5} Toms, A.\ S.: {\it Stability in the Cuntz semigroup}, arXiv preprint math.OA/0607099 (2006)

\bibitem{To6} Toms, A.\ S.: {\it Flat dimension growth for C$^*$-algebras}, J. Funct. Anal. {\bf 238} (2006), 678-708

\bibitem{To7} Toms, A.\ S.: {\it An infinite family of non-isomorphic C$^*$-algebras with identical $\mathrm{K}$-theory}, arXiv preprint
math.OA/0609214 (2006)

\bibitem{TW1} Toms, A.\ S.\ and Winter, W.: {\it Strongly self-absorbing $C^{*}$-algebras}, 
arXiv preprint math.OA/0502211 (2005); Trans.\ Amer.\ Math.\ Soc., to appear

\bibitem{TW2} Toms, A.\ S.\ and Winter, W.: {\it $\mathcal{Z}$-stable ASH algebras}, 
arXiv preprint  math.OA/0508218 (2005); Canad.\ J.\ Math., to appear

\bibitem{V1} Villadsen, J.: {\it Simple $C^*$-algebras with perforation},
J.\ Funct.\ Anal.\ {\bf 154} (1998), 110-116

\bibitem{V2} Villadsen, J.: {\it On the stable rank of simple 
$C^{*}$-algebras}, J.\ Amer.\ Math.\ Soc.\ {\bf 12} (1999), 1091-1102

\bibitem{Wi1} Winter, W.: {\it Covering dimension for nuclear 
$C^*$-algebras},  J.\ Funct.\ Anal.\ {\bf 199} (2003), 535-556

%\bibitem{Wi2} Winter, W.: {\it Covering dimension for nuclear $C^*$-algebras II}, preprint math.OA/0108102, Trans.\ Amer.\ Math.\ Soc., to appear 

%\bibitem{Wi3} Winter, W.: {\it Decomposition rank of subhomogeneous 
%$C^*$-algebras},  Proc.\ London 
%Math.\ Soc.\ {\bf 89} (2004), 427-456

\bibitem{Wi4} Winter, W.: {\it On topologically finite-dimensional simple 
$C^{*}$-algebras}, Math.\ Ann.\ {\bf 332} (2005), 843-878

\bibitem{Wi5} Winter, W.: {\it On the classification of simple $\mathcal{Z}$-stable $C^{*}$-algebras with real rank 
zero and finite decomposition rank}, J.\ London Math.\ Soc.\ {\bf 74} (2006), 167-183

\bibitem{Wi6} Winter, W.: {\it Simple $C^{*}$-algebras with locally finite decomposition rank}, arXiv preprint math.OA/0602617 (2006)




\end{thebibliography}
\end{document}